\documentclass[final,5p,times]{elsarticle}

\biboptions{sort&compress,numbers}
\usepackage[nodots]{numcompress}
\usepackage{lineno,url}
\usepackage{amsmath,amsfonts,amssymb}

\usepackage{amssymb,amsmath,amsthm}
\usepackage{overpic}
\usepackage{contour}\contourlength{0.7pt}
\usepackage{tikz}\usetikzlibrary{shapes,arrows}

\usepackage{versions}
\includeversion{PDF}

\newtheoremstyle{mytheorem}{5pt plus 5pt minus 3pt}{4pt plus 3pt minus 1.5pt}
	{\itshape}{}{\bfseries}{.}{1ex plus 1ex minus .5ex}{}
\newtheoremstyle{mydef}{5pt plus 5pt minus 3pt}{4pt plus 3pt minus 1.5pt}
	{}{0pt}{\bfseries}{.}{1ex plus 1ex minus .5ex}{}
\newtheoremstyle{myremark}{5pt plus 5pt minus 3pt}{4pt plus 3pt minus 1.5pt}
	{}{0pt}{\itshape}{.}{1ex plus 1ex minus .5ex}{}
\theoremstyle{mytheorem}
\newtheorem{prop}{Proposition}[section]
\newtheorem{lemma}{Lemma}[section]
\theoremstyle{mydef}
\newtheorem{dfn}[prop]{Definition}
\theoremstyle{myremark}

\newtheorem{rem}{Remark}

\newcommand{\bc}{\mathbf c}
\newcommand{\bw}{\mathbf w}
\newcommand{\bm}{\mathbf m}
\newcommand{\bn}{\mathbf n}
\newcommand{\bd}{\mathbf d}
\newcommand{\be}{\mathbf e}
\newcommand{\bp}{\mathbf p}

\newcommand{\bv}{\mathbf v}
\newcommand{\bu}{\mathbf u}
\newcommand{\bs}{\mathbf S}

\newcommand{\mc}{\mathcal C}
\newcommand{\bj}{\mathbf J}
\newcommand{\bk}{\mathbf K}
\newcommand{\bq}{\mathbf q}
\newcommand{\bQ}{\mathbf Q}

\def\R {\mathbb R}

\newcommand{\dpr}[2]{{\langle #1, #2 \rangle}}

\def\Acaption#1#2{\caption{#2}\vspace*{-#1}}

\definecolor{Mgreen}{RGB}{34,139,34}
\definecolor{blau}{rgb}{0.15,0.2,0.5}
\definecolor{gray}{rgb}{0.5,0.5,0.5}
\definecolor{drot}{rgb}{0.7,0,0.1}
\definecolor{gelb}{rgb}{.55,.40,.1}

\begin{document}

\sloppy
\begin{frontmatter}

\title{Strain-minimising Stream Surfaces}

\author[Kaust]{Michael Barto\v{n}}
\ead{Michael.Barton@kaust.edu.sa}
\author[Cam]{Ji\v{r}\'{i} Kosinka}
\ead{Jiri.Kosinka@cl.cam.ac.uk}
\author[Kaust]{Victor M. Calo}
\ead{Victor.Calo@kaust.edu.sa}

\address[Kaust]{Numerical Porous Media Center, King Abdullah University of Science and Technology,
                Thuwal 23955-6900, Kingdom of Saudi Arabia}
\address[Cam]{Computer Laboratory, University of Cambridge,
              15 JJ Thomson Avenue, Cambridge CB3 0FD, United Kingdom}

\begin{abstract}
We study the problem of finding strain-minimising stream surfaces in a
divergence-free vector field. These surfaces are generated by motions
of seed curves that propagate through the field in a strain
minimising manner, i.e., they move without stretching or shrinking,
preserving the length of their arbitrary arc. In general fields, such curves do not
exist. However, the divergence-free constraint gives rise to
these 'strain-free' curves that are locally arc-length preserving when
infinitesimally propagated. Several families of strain-free curves are
identified and used as initial guesses for stream surface generation.
These surfaces are subsequently globally optimised to obtain the best
strain-minimising stream surfaces in a given divergence-free vector
field.

Our algorithm was tested on benchmark datasets, proving its
applicability to incompressible fluid flow simulations, where our
strain-minimising stream surfaces realistically reflect the flow of a
flexible univariate object.
\end{abstract}


\end{frontmatter}

\section{Introduction}\label{sec:intro}

We investigate a special class of stream surfaces generated
by seed curves that minimise a certain arc-length energy.
With flow simulations of viscous fluids in mind, consider a seed curve as a sequence of infinitely many liquid drops.
We seek to answer the following question: Given a time-independent vector field in 3D, is there a curve that propagates in time without changing the strain between any two of its neighbouring drops? In other words, does a curve
that moves in the vector field without shrinking or stretching any of its arcs exist?

It is well known
\cite{Davis-1967-VF}
that divergence-free vector fields are volume preserving, i.e., the volume of any 3D object remains constant when propagated in time. In general, this is not true for lower-dimensional objects. Thus,  a natural question arises: Are there lower-dimensional objects (surfaces, curves) that propagate in time in the same manner, i.e., preserving their lower-dimensional measures (area, length)? Or, stated differently, is there a region where the divergence-free vector field acts on an object more than in volume preserving manner, namely by preserving its 1D or 2D measure? This paper investigates this question in the curve case.

Our motivation is straightforward: if a sufficiently elastic univariate object (an elastic rod or drops of another fluid) is put into the flow along a special curve, the deformation that acts on it is bending only; it contains no `strain' forces.

\paragraph*{Problem formulation.}
Given a steady (time independent), divergence-free vector field $\bv$ in some domain $\Omega \subset \R^3$, ($\nabla\cdot \bv(\bp) = 0, \forall \bp \in \Omega $),
find a \emph{stream surface} such that the \emph{seed curve}
that defines it propagates in time along the surface so that its arc-length changes as little as possible.

We combine theoretical investigations and a practical algorithm for finding such stream surfaces. The main steps and contributions of our method are:
\begin{itemize}
\item We theoretically investigate families of seed curves based on certain strain-minimising energies (Section~\ref{sec:SE}).
\item These candidate seed curves are used to generate initial stream surfaces (Section~\ref{sec:smss}).
\item The initial stream surfaces are globally optimised and ranked according to their strain energies (Section~\ref{sec:opt}).
\end{itemize}
Our implementation of the method is presented in Section~\ref{sec:impl}.
We have validated our theoretical results on several benchmark datasets and demonstrated the applicability of our method on numerous examples (Section~\ref{sec:exmp} and the accompanying video).
Possible extensions of our method are discussed in Section~\ref{sec:dis} and the paper is concluded in Section~\ref{sec:conclu}.

\section{Related work}\label{sec:rel}

Stream surfaces, used as a tool for visualising characteristic features of vector fields, have been extensively studied in the visualisation literature; see
\cite{Peikert-2002-Turbine,Theisel-2013-StreamSrf,Edmunds-2012-AutomStrSrf,Reviewer3Citation}, the survey paper \cite{McLoughlin-2010-Survey} and the references cited therein.
Classical methods \cite{Hultquist-1992-SeedingLine} are usually based on trial-and-error approaches: the user inserts seed curves (typically straight lines), stream surfaces are computed, and, if they do not capture desired features well, the initial seed curves are modified and the whole process is repeated. Since visualising vector fields by stream surfaces (compared to using streamlines) has became more popular \cite{McLoughlin-2010-Survey},
research in \emph{automatic} stream surface seeding has recently become very active. Our method fits in this modern family of automatic stream surface algorithms.

Divergence-free vector fields are used in many areas and applications such as incompressible fluid simulations \cite{Fedkiw-2010-FluidSimul}, smoke visualisation \cite{Weinkauf-2008-SquareCylinder}, and are also a favourite modelling/deformation tool \cite{Seidel-2006-Deformations} due to their volume preserving property.

Another research area, which our work is connected to, relates to curve evolution \cite{Younes-2010}, where typically the curve and the property
to preserve (e.g. arc-length) are given and the corresponding evolving vector fields are sought after. For example, \cite{Polthier-2005-CrvEvolution} seeks
3D smoothing flows that satisfy additional spatial constraints, and \cite{Barton-2012-Snakes} uses the flow to evolve curves in 3D space while preserving their arc-length and curvature.

\section{Strain-energy minimising curves}\label{sec:SE}

\begin{figure}
\begin{PDF}
 \centerline{\hfill
 \begin{overpic}[width=.69\columnwidth,angle=0]{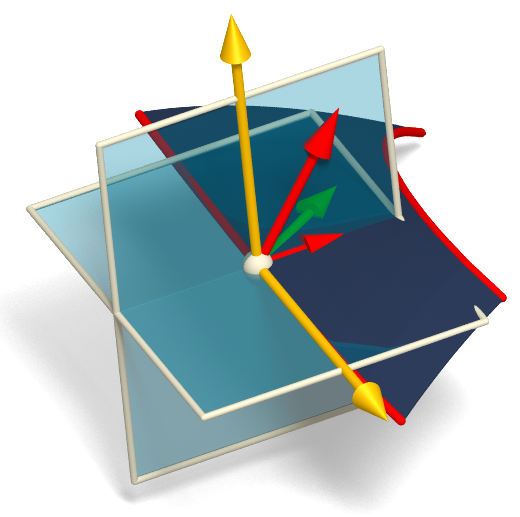}
    \put(65,80){$\bs_{st}$}
    \put(65,13){$\bs_{s}=\bu$}
    \put(50,90){$\bn$}
    \put(42,42){$\bp$}
    \put(5,63){$\tau_{\bp}$}
    \put(90,30){$\bs(s,t)$}
    \put(66,61){\color{Mgreen}$\bv$}
    \put(10,80){$\gamma(s) = \bs(s,0)$}
	\end{overpic}
 \hfill}
 \vspace{-4pt}
\end{PDF}
 \Acaption{1ex}{Strain minimisation in terms of differential geometry.
 At any point $\bp$ on the seed curve $\gamma(s) = \bs(s,0)$, the directional derivative $\bs_{st}$ of the vector $\bu$ along $\bv$
is perpendicular to $\bu$; see \eqref{eq:conjugate}.}\label{fig:DiffGeo}
 \end{figure}

Let $\bv(\bp)$ be a steady differentiable vector field defined over a domain $\Omega\subset\R^3$ and assume that $\bv$ is divergence free, i.e., $\mathrm{div}\,\, \bv(\bp) = \dpr{\nabla}{\bv(\bp)} = 0, \forall \bp \in \Omega$.
Let $\bj$ be the Jacobian matrix of $\bv(\bp)$, i.e., $\bj_{ij}=\frac{\partial \bv_i}{\partial p_j}$ with $\bp=(p_1,p_2,p_3)$.

Consider a regular curve $\gamma(s)$ parametrised by arc-length, $s\in[s_0,s_1]$. We regard $\gamma$ as a \emph{seed curve} that gives rise to a \emph{stream surface} $\bs(s,t)$ with normal $\bn(s,t)$, i.e., $\dpr{\bv(\bs(s,t))}{\bn(s,t)}=0$ for all $(s,t)$ in the surface domain $[s_0,s_1]\times[t_0,t_1]$,
and $\bs(s,0)=\gamma(s)$, $0\in[t_0,t_1]$; see Fig.~\ref{fig:DiffGeo}. The partial derivatives of $\bs$ will be denoted $\bs_{s}$, $\bs_{st}$, etc.

In general, curves which maintain their arc-length constant (i.e., equal to $s_1-s_0$ independently of $t$) when deformed by $\bv$ do not exist. However, as we show below, it is possible to find curves which approximate this property to first or even second order. To make this concept precise, we formulate the following
\begin{lemma}\label{le:Taylor}
The Taylor expansion of the arc length of $\gamma$ with respect to $\bv$ at $t=0$ is given by
\begin{equation}
\int_{s_0}^{s_1} ||\bs_s(s,t)|| \, \mathrm{d}s = (s_1-s_0) + c_1 t + c_2 t^2 + \mathcal{O}(t^3)
\end{equation}
with
\begin{equation}\label{eq:c12}
\begin{array}{rcl}
c_1 & = & \int_{s_0}^{s_1}\dpr{\bs_s}{\bs_{st}} \, \mathrm{d}s |_{t=0},\\
c_2 & = & \int_{s_0}^{s_1}\dpr{\bs_{st}}{\bs_{st}}+\dpr{\bs_s}{\bs_{stt}}-\dpr{\bs_s}{\bs_{st}}^2 \, \mathrm{d}s |_{t=0}.
\end{array}
\end{equation}
\end{lemma}
A straightforward proof can be found in Appendix A.

\begin{figure}
\begin{PDF}
 \centerline{\hfill
 \begin{overpic}[width=.89\columnwidth]{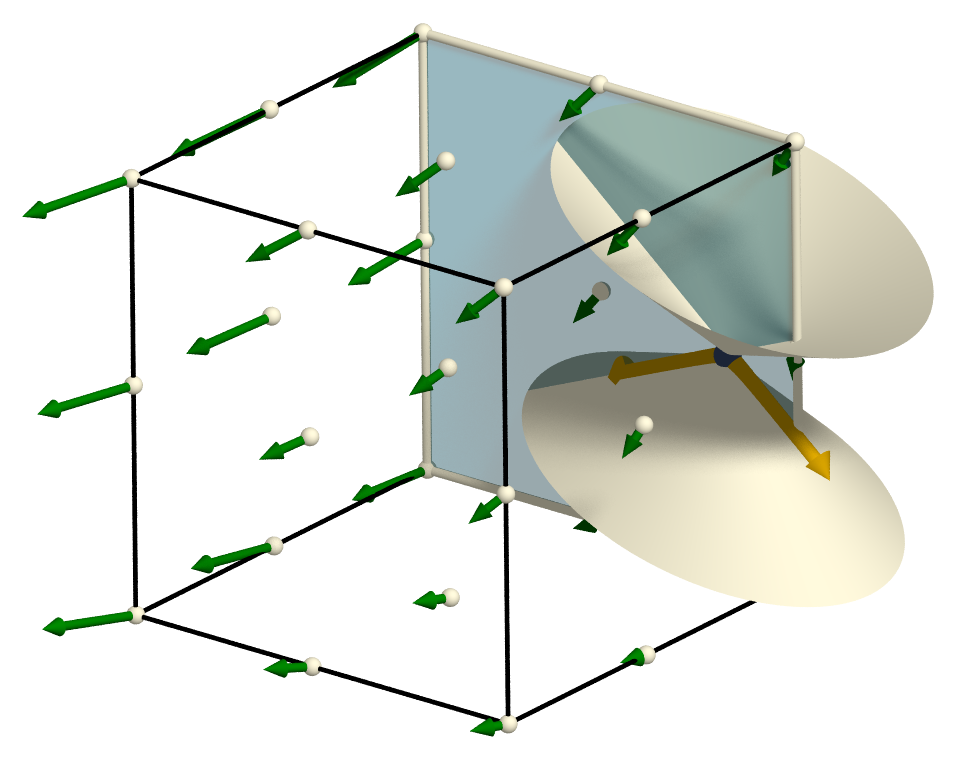}
    \put(85,25){$\bd_1$}
    \put(62,37){$\bd_2$}
    \put(72,38){$\bp$}
    \put(70,70){$\partial \Omega$}
    \put(80,10){$\mc$}
	\end{overpic}
 \hfill
 \vspace{-6pt}}
\end{PDF}
 \Acaption{1ex}{Computing first-order vectors on the boundary of $\Omega$. Self-conjugate vectors associated with $\bj$ (see
 \eqref{eq:conjugate}) form a quadratic cone $\mc$. If restricted to the tangent space of $\Omega$ at ${\bp}$, at most two first-order vectors $\bd_1$ and $\bd_2$ exist.
}\label{fig:Cone}
 \end{figure}

Our aim is to identify curve(s) $\gamma$ in $\Omega$ for which $c_1$ (and also $c_2$, if possible) vanishes. These will subsequently be used to identify strain-minimising stream surfaces in an optimisation procedure.

\subsection{First-order strain energy}\label{sec:SE1}

\begin{figure*}
\begin{PDF}
  \hfill
 \begin{overpic}[width=.33\textwidth]
	{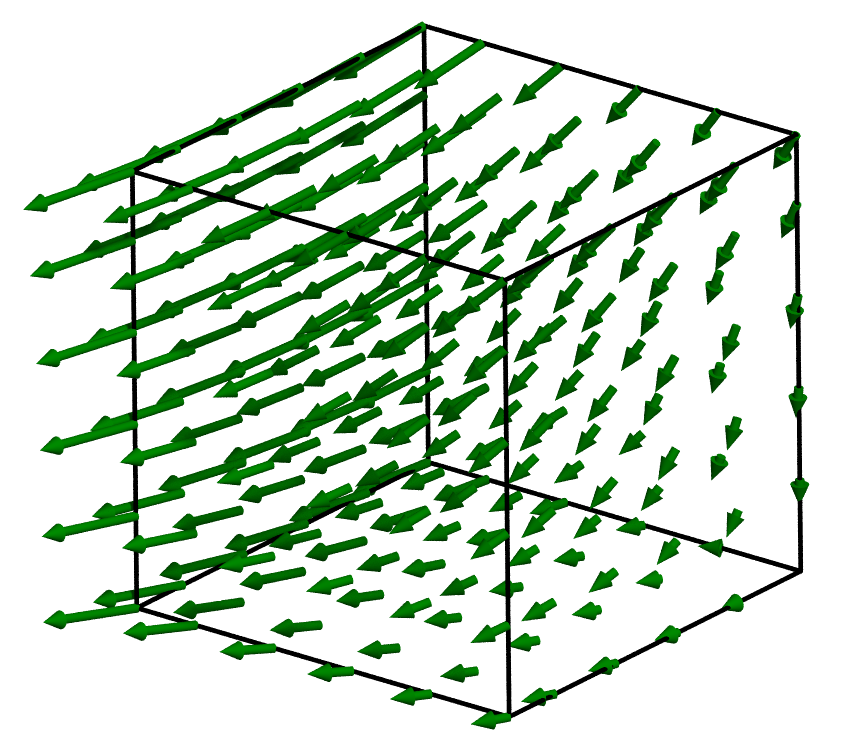}
	\put(0,0){(a)}
	\end{overpic}\hfill
  \begin{overpic}[width=.33\textwidth]
	{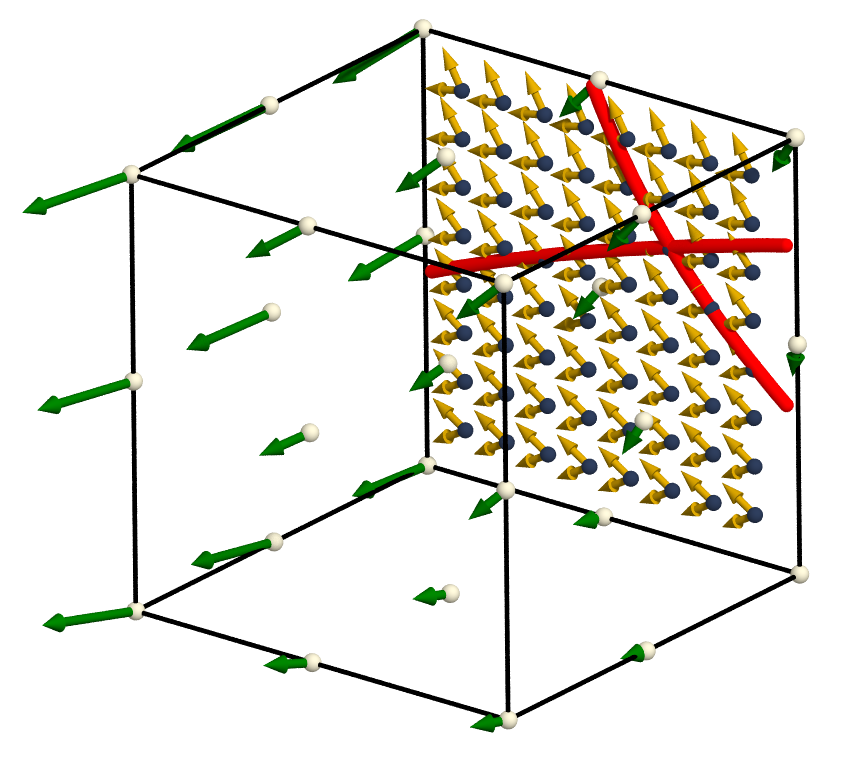}
	\put(0,0){(b)}
    \put(72,82){$\gamma_1$}
    \put(95,60){$\gamma_2$}
	\end{overpic}\hfill
 \begin{overpic}[width=.33\textwidth]
	{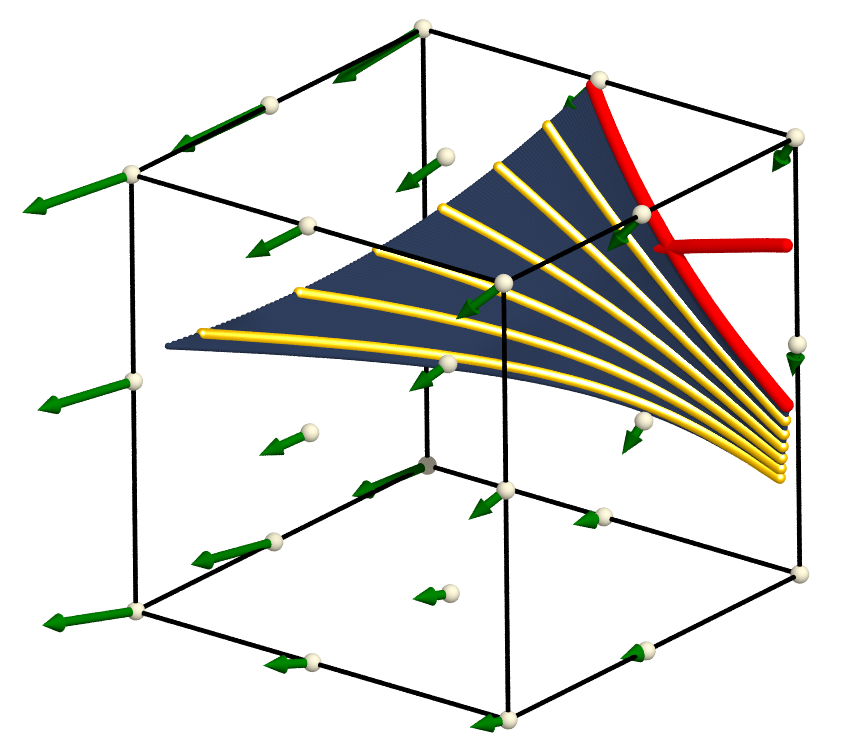}
	\put(0,0){(c)}
	\end{overpic}
  \hfill
  \vspace{-3pt}
\end{PDF}
  \Acaption{1.5em}{Initialisation of first-order curves. (a) A divergence-free vector field $\bv=[x+y^2+2z^3,10x^3+2y,2x^2y-3z]$
within a unit domain $\Omega$. (b) Two first-order vector fields (yellow) on the boundary quad given by $y=0$. Integration gives rise to first-order boundary
curves ($\gamma_1$ and $\gamma_2$ in red). (c) The stream surface obtained by propagating $\gamma_1$ in $\bv$.
Observe that the preservation of arc length is guaranteed
only locally, i.e., close to the boundary, up to first order via~\eqref{eq:conjugate}.}
\label{fig:Initialization}
  \end{figure*}

We have $\bs_{t} = \bv$ and we denote $\bu:=\dot{\gamma}(s) = \bs_{s}(s,0)$.
Then
\begin{equation}\label{eq:1stOrderAnal}
 \dpr{\bs_{s}}{\bs_{s}} |_{t=0} \equiv \dpr{\bu}{\bu} \equiv 1.
\end{equation}
Differentiation with respect to $t$ yields (see Fig.~\ref{fig:DiffGeo})
\begin{equation}\label{eq:conjugate}
 \dpr{\bs_{s}}{\bs_{st}} |_{t=0} \equiv \bu \bj \bu^\top \equiv 0.
\end{equation}
This observation leads us to

\begin{dfn}\label{def:strainE}
The \emph{first-order strain energy} of $\gamma$ is given by
\begin{equation}\label{eq:E1}
E_1(\gamma) = \frac{1}{s_1-s_0}\int_{s_0}^{s_1} \dpr{\bs_{s}}{\bs_{st}}^2 \, \mathrm{d}s |_{t=0}.
\end{equation}
A curve $\gamma$ on which $E_1$ vanishes will be called a \emph{first-order curve}, and their collection denoted
$$
\Gamma_1=\{\gamma\, | \, \gamma\subset\Omega, E_1(\gamma) = 0\}.
$$
Moreover, a vector $\bd$ such that $\bd \bj \bd^\top=0$, i.e., a self-conjugate vector associated with $\bj$, will be called a \emph{first-order vector}.
\end{dfn}

Note that $E_1(\gamma) = 0$ implies $c_1=0$ and thus the first order energy is well defined. In other words, $E_1$ measures the change of the magnitude of $\dot{\gamma}$ along $\gamma$ for an infinitesimal increment of $t$. $E_1(\gamma) = 0$ implies the deformation given by $\bv$ that acts on $\gamma$ preserves, up to \emph{first order}, the magnitude of its tangent vector $\dot{\gamma}$ and hence is locally arc-length preserving.

In order to find seed curves that belong to $\Gamma_1$ with respect to $\bv$, we seek first-order vectors. To this end, we denote \emph{strain rate} $\bj^+ := (\bj+\bj^\top)/2$ and \emph{vorticity} $\bj^- := (\bj-\bj^\top)/2$ as the symmetric and antisymmetric parts of the Jacobian matrix of $\bv$. A vector field $\bv$ for which $\bj^+$ is regular will be called \emph{non-degenerate}.

\begin{lemma}\label{lem:dirs}
Let $\bp\in\Omega$ and $\bv$ be divergence free and non-degenerate in a neighbourhood of $\bp$. Then all first-order vectors $\bd$ form a quadratic cone with apex at $\bp$.
\end{lemma}
\noindent
\textbf{Proof.}
Since $\bu \bj^- \bu^\top=0$ for any $\bu$, we obtain the condition $\bd \bj^+ \bd^\top=0$. By definition, $\mathrm{div}\,\, \bv(\bp)=0$ is equivalent to $tr(\bj) = 0$ at $\bp$ and thus $tr(\bj^+) = 0$. From the spectral theorem it follows that $\bj^+$ has three real eigenvalues such that $\lambda_1+\lambda_2+\lambda_3 = 0$. This in turn implies that the signature of $\bj^+$ is either $(+,+,-)$ or $(+,-,-)$ by non-degeneracy of $\bv$. Consequently, all solutions of $\bd \bj \bd^\top=0$ form a quadratic cone.
\hfill $\square$

\begin{rem}\label{rem:singular}
In the special case (of measure zero) when $\bj^+$ is singular, the space of first-order vectors $\bd$ at $\bp$ that solve $\bd\bj\bd^\top=0$ is either given by two planes intersecting in a line incident with $\bp$ (the signature of $\bj^+$ is $(+,-,0)$) or any vector is a first-order vector (the signature of $\bj^+$ is $(0,0,0)$). Consequently, first-order vectors exist at any point $\bp\in\Omega$ for divergence-free vector fields.
\end{rem}

From Lemma~\ref{lem:dirs} and Remark~\ref{rem:singular} it follows that there exist infinitely many first-order curves $\gamma$ passing through every point in $\Omega$. They can be obtained by integrating first-order vectors $\bd$, which, however, form a multi-valued field. Thus, the set $\Gamma_1$ of these curves is too large to be practical. We therefore explore three conditions that select special classes of first-order curves from $\Gamma_1$:
\begin{enumerate}
 \item minimise a certain second-order strain energy;
 \item restrict $\Gamma_1$ to curves on the boundary only, i.e., $\gamma\subset\partial\Omega$;
 \item constrain the variation of $\bd=\dot{\gamma}$ along $\gamma$.
\end{enumerate}
We now address each of these strategies in detail.


\subsection{Second-order strain energy}\label{sec:SE2}

\begin{figure}
\begin{PDF}
  \hfill
 \begin{overpic}[width=.22\textwidth]
	{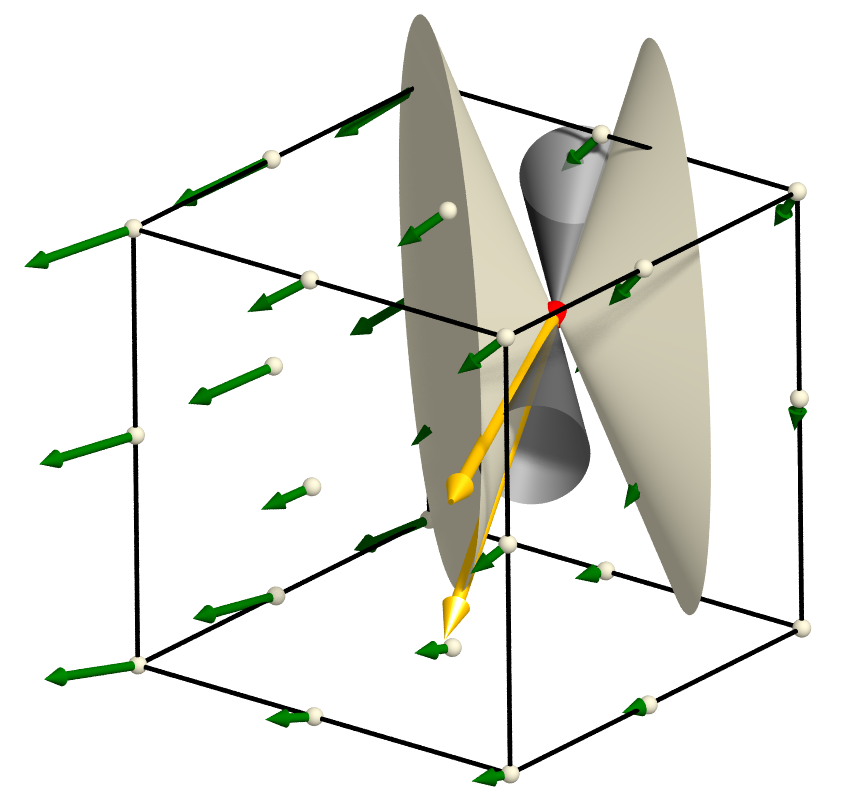}
    \put(20,80){$\Omega$}
	\end{overpic}\hfill
  \begin{overpic}[width=.22\textwidth]
	{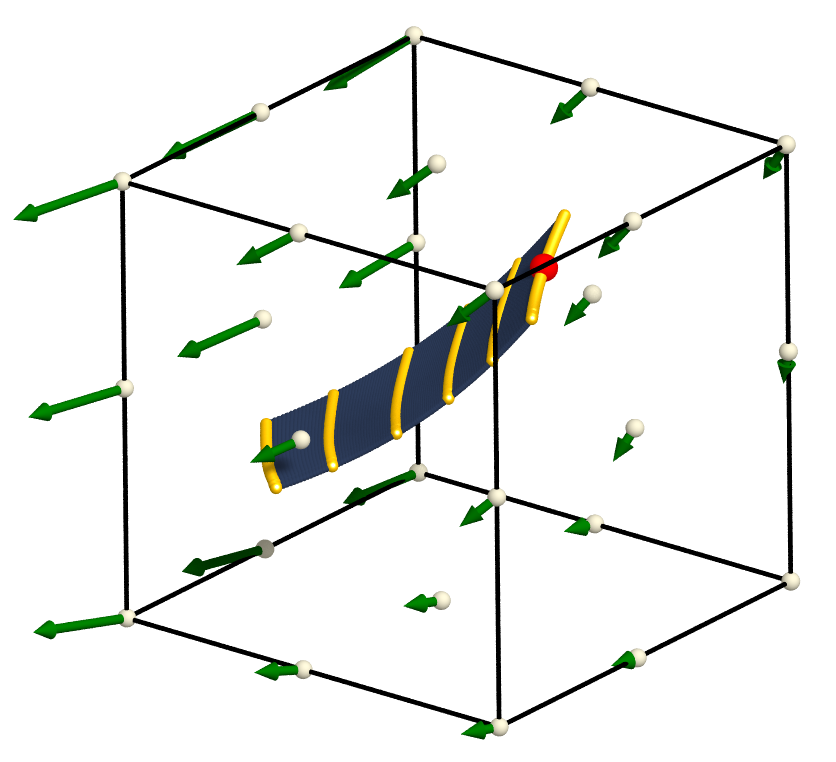}
    \put(65,45){$\gamma$}
	\end{overpic}\hfill
\end{PDF}
  \Acaption{1.5em}{Second-order curves. Left: At an interior point (red) of $\Omega$, second-order vectors, generically and if they exist, correspond to the intersection of two quadratic cones. Right: Integrating second-order vectors gives $\gamma$, an integral curve whose tangent vectors $\dot{\gamma}$
solve \eqref{eq:conjugate} and \eqref{eq:2ndOrderAnal}. The stream surface emanating from $\gamma$ for the same vector field as in
Fig.~\ref{fig:Initialization} is shown.}\label{fig:Cone2}
  \end{figure}

\begin{figure*}
\begin{PDF}
  \hfill
 \begin{overpic}[width=.33\textwidth]
	{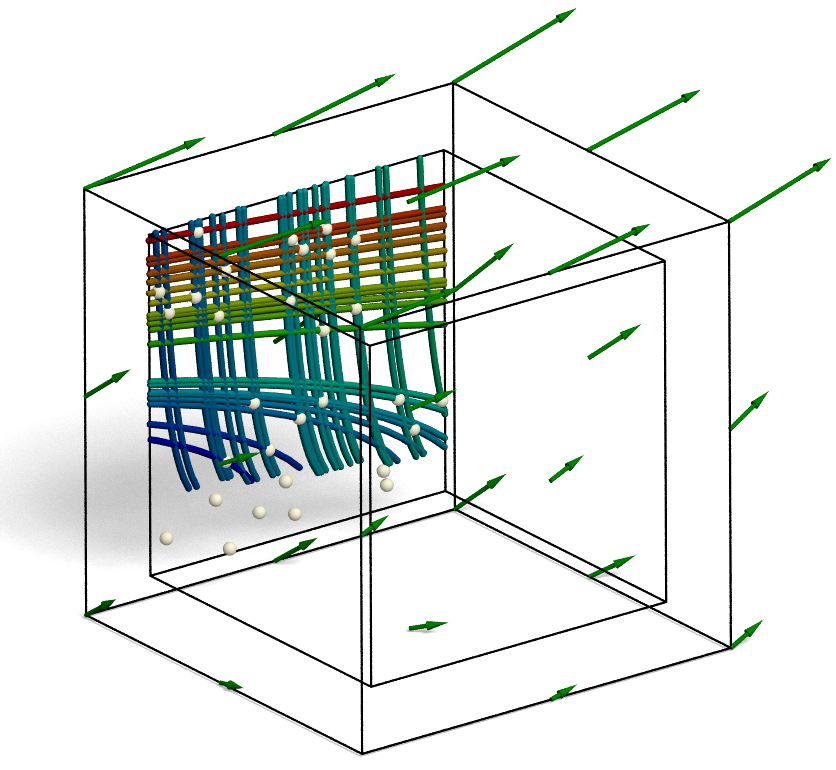}
    \put(0,80){(a)}
    \put(40,72){$\partial \Omega$}
    \put(30,80){$\partial\overline{\Omega}$}
	\end{overpic}\hfill
  \begin{overpic}[width=.33\textwidth]
	{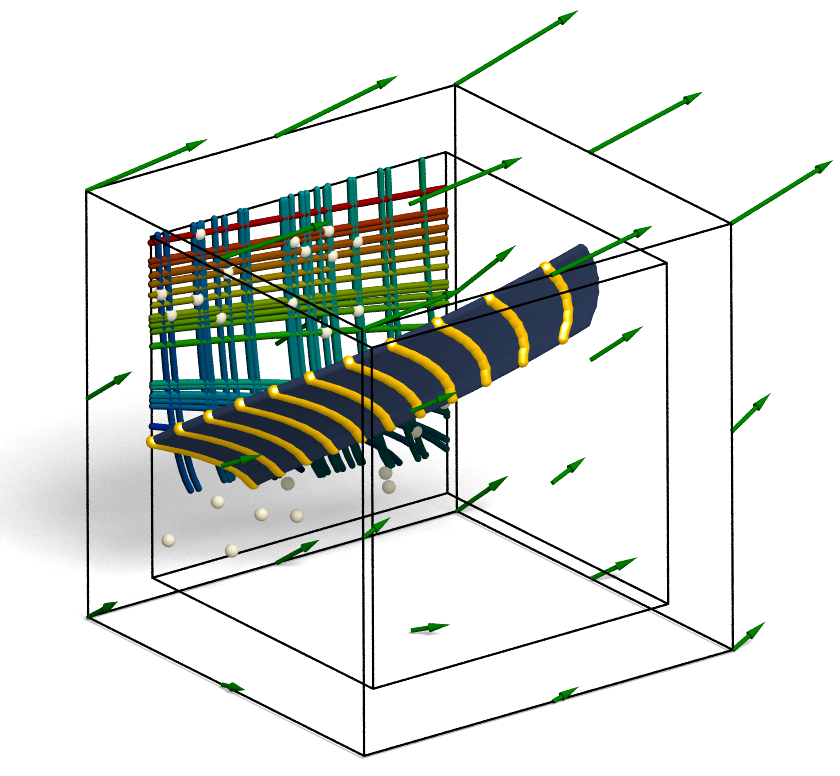}
    \put(-17,-1){\includegraphics[width=0.1\textwidth]{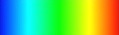}}
    \footnotesize
	\put(-22,-6){$\min E_2$}
    \put(8,-6){$\max E_2$}
	\normalsize
	\put(0,80){(b)}
	\end{overpic}\hfill
 \begin{overpic}[width=.33\textwidth]
	{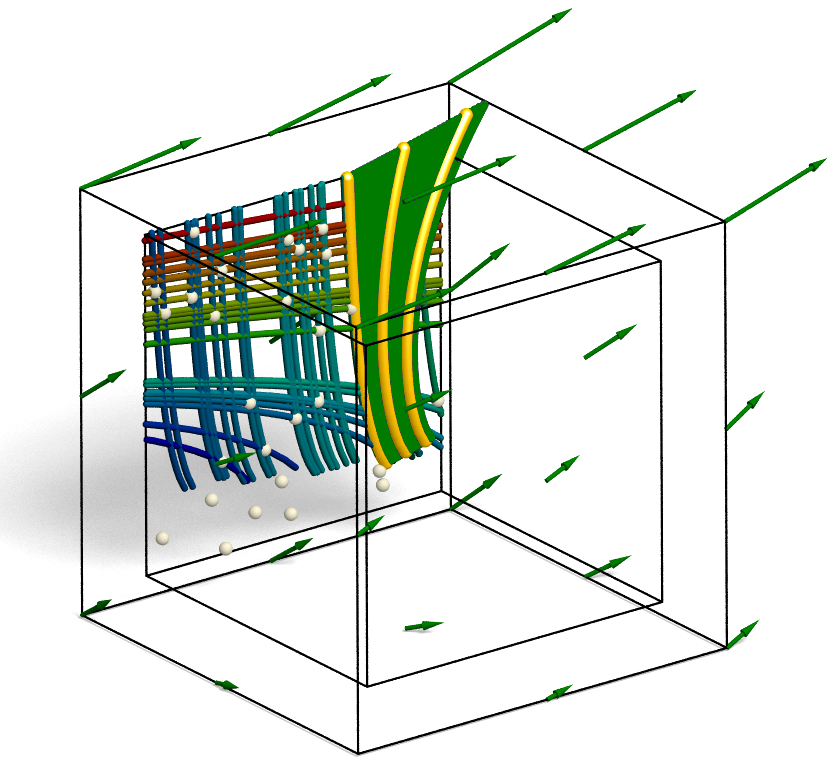}
	\put(0,80){(c)}
	\end{overpic}\hfill
\vspace{-0.01cm}
 \hfill
 \begin{overpic}[width=.24\textwidth]
	{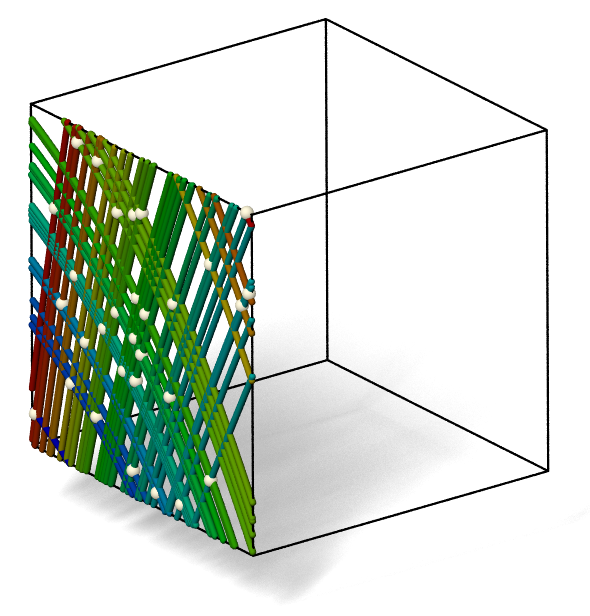}
	\end{overpic}\hfill
  \begin{overpic}[width=.24\textwidth]
	{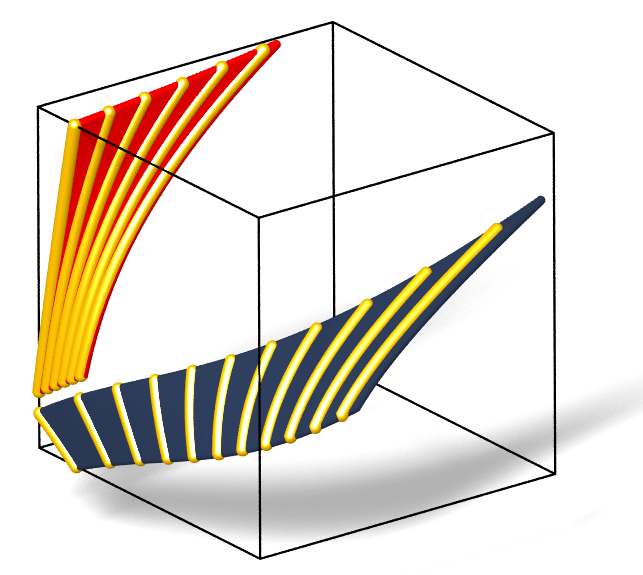}
	\put(-10,2){(d)}
	\end{overpic}\hfill
 \begin{overpic}[width=.24\textwidth]
	{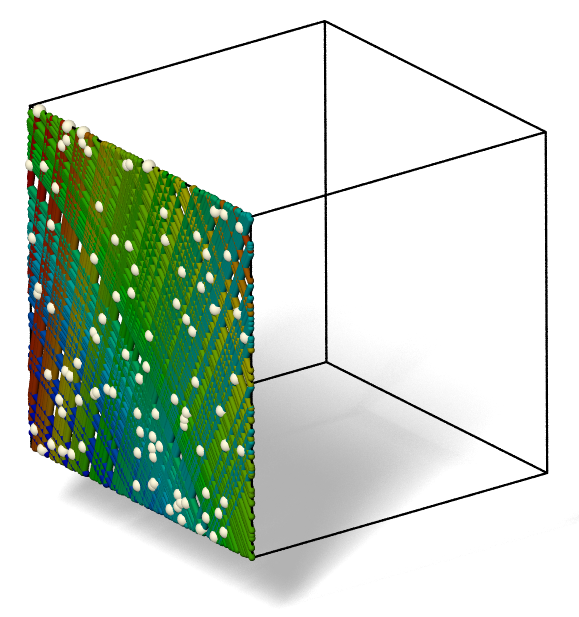}
	\end{overpic}\hfill
\begin{overpic}[width=.24\textwidth]
	{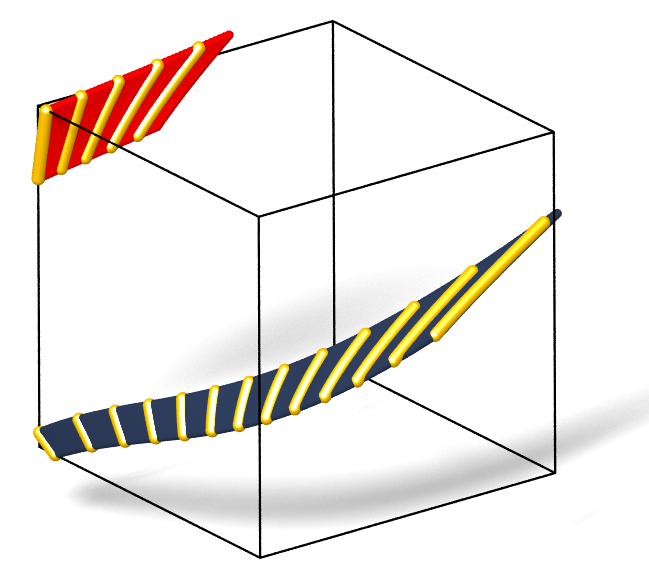}
	\put(-10,2){(e)}
	\end{overpic}\hfill
    \vspace{-3pt}
\end{PDF}
  \Acaption{1.5em}{Exploring the space of first-order boundary curves.
	(a) At sampled boundary points, the curves are computed and colour-coded according to the second-order energy $E_2$; see \eqref{eq:E2}.
	(b)--(c) The stream surface obtained from the curves corresponding to the best and the median value of $E_2$, respectively. Seed curve integration is terminated when the outer offset of $\partial \Omega$, $\partial\overline{\Omega}$, is reached.
	(d)--(e) The analogous situation for a different part of $\partial\Omega$. Left: First-order boundary curves arisen from $n=30$ and $100$ random samples, respectively. Right: The stream surfaces obtained by integrating the vector field starting from the best (blue) and the worst (red) first-order boundary curves according to $E_2$.}\label{fig:EnergyRanking}
  \end{figure*}

Similarly to the first-order strain energy $E_1$, we now define its second-order counterpart. Differentiating \eqref{eq:conjugate} with respect to $t$ gives
\begin{equation}\label{eq:2ndOrderAnal}
 \langle \bs_{st},\bs_{st} \rangle |_{t=0} +
 \langle \bs_{s},\bs_{stt} \rangle |_{t=0} \equiv
 \bu\bk\bu^\top \equiv 0,
\end{equation}
where $\bk = \bj\bj^\top+\bj^2$. This leads us to
\begin{dfn}\label{def:strainE2}
The \emph{second-order strain energy} of $\gamma$ is given by
\begin{equation}\label{eq:E2}
E_2(\gamma) = \frac{1}{s_1-s_0}\int_{s_0}^{s_1}
(\langle \bs_{st},\bs_{st} \rangle + \langle \bs_{stt}, \bs_{s} \rangle)^2 \, \mathrm{d}s |_{t=0}.
\end{equation}
A curve $\gamma$ on which both $E_1$ and $E_2$ vanish will be called a \emph{second-order curve}, and their collection denoted
$$
\Gamma_2=\{\gamma\, | \, \gamma\subset\Omega, E_1(\gamma) = E_2(\gamma) = 0\}.
$$
Moreover, a first-order vector $\bd$ such that $\bd\bk\bd^\top = 0$ will be called a \emph{second-order vector}.
\end{dfn}

Note that $E_1(\gamma) = E_2(\gamma)= 0$ implies that $c_1=c_2=0$ in Lemma~\ref{le:Taylor} and thus the second-order energy is well defined: second-order curves are characterised by the property that their deformation given by $\bv$ that acts on $\gamma$ preserves the magnitude of the tangent vector up to \emph{second order} at $t=0$; cf. Lemma~\ref{le:Taylor}.

Whereas the existence of first-order vectors $\bd$ is guaranteed at every $\bp$
(see Fig.~\ref{fig:Cone} and Lemma~\ref{lem:dirs}), a non-zero solution of \eqref{eq:2ndOrderAnal} need not exist.
The generic cases with respect to $\bd\bk\bd^\top=0$ are categorised by the signature of $\bk^+$: $(+,+,+)$ or $(-,-,-)$ yield no non-zero solution; $(+,+,-)$ or $(+,-,-)$ give, as in the first-order case, a quadratic cone. In the former case, there are no second-order vectors at $\bp$ and alternatives must be sought; see Sections~\ref{sec:bnd} and \ref{sec:str}. The latter case leads to the intersection of two quadratic cones. The situation is shown in Fig.~\ref{fig:Cone2}, left. Up to four second-order vectors can be found by solving a quartic equation or, more geometrically, by reducing the problem via a cubic equation to the intersection of a quadratic cone with two planes; see Appendix B. Non-generic cases are, for the sake of brevity, not considered.

Therefore, testing whether a real non-zero solution of the system given by \eqref{eq:conjugate} and \eqref{eq:2ndOrderAnal} exists is a cheap closed-form operation. This allows us to quickly explore $\Omega$ for regions where second-order vectors exist and, in the positive case, to integrate them to obtain curves in $\Gamma_2$. An example of a second-order curve is shown in Fig.~\ref{fig:Cone2}, right. Such curves then form input for our optimisation algorithm (Section~\ref{sec:opt}).

If no second-order curves exist in $\Omega$ or if the existing ones are not satisfactory for visualisation purposes, we employ our second strategy: boundary curves.

\subsection{First-order boundary curves}\label{sec:bnd}

Our second strategy is to restrict the set $\Gamma_1$ to curves on the boundary of $\Omega$. This is a reasonable restriction from the point of view of visualisation: boundary seed curves capture the behaviour of a given flow as it enters/exits the domain; see Fig.~\ref{fig:Initialization}. When the point $\bp$ lies on the boundary $\partial\Omega$ of $\Omega$ and $\bv$ is non-degenerate, we see that there exist at most two first-order vectors $\bd_1$ and $\bd_2$ in the tangent space
of $\partial \Omega$
at $\bp$; see Fig.~\ref{fig:Cone}.

\emph{First-order boundary curves} form a subset of $\Gamma_1$, which we denote
\begin{equation}\label{eq:family1}
\Gamma_1^\partial = \{\gamma \, | \, \gamma \subset \partial \Omega, E_1(\gamma)  = 0 \}.
\end{equation}
The first-order strain energy $E_1$ was defined to measure the change of the magnitude of the unit tangent vector $\dot{\gamma}$ of $\gamma$,
when being instantaneously moved by the vector field $\bv$.
Every curve in $\Gamma_1^\partial$ is a first-order curve,
but is not, in general, a second-order curve.
Nevertheless, $E_2$ can be used as a ranking criterion to determine
good candidate seed curves among those in $\Gamma_1^\partial$; see Fig.~\ref{fig:EnergyRanking}.

\begin{figure*}[!tbh]
\begin{PDF}
\begin{overpic}[width=.25\textwidth]
    {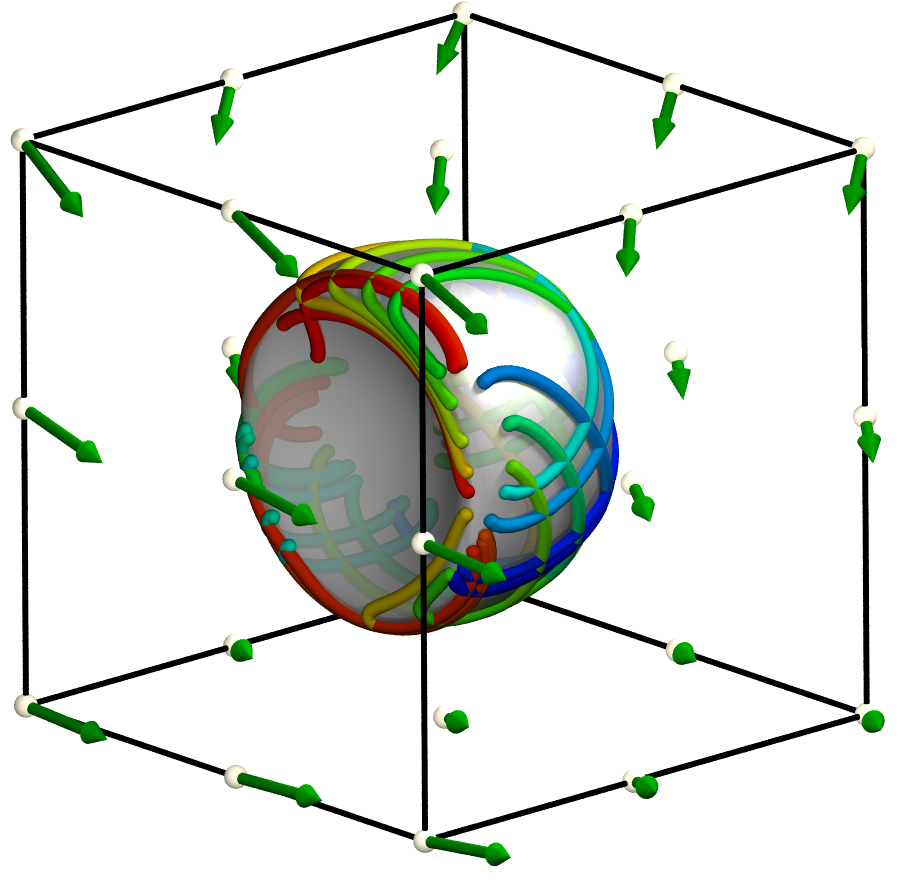}
    \put(50,21){$\partial\Omega$}
    \put(0,-5){(a)}
\end{overpic}\hfill
\raisebox{0.0\columnwidth}{
\begin{overpic}[width=.2\textwidth]
    {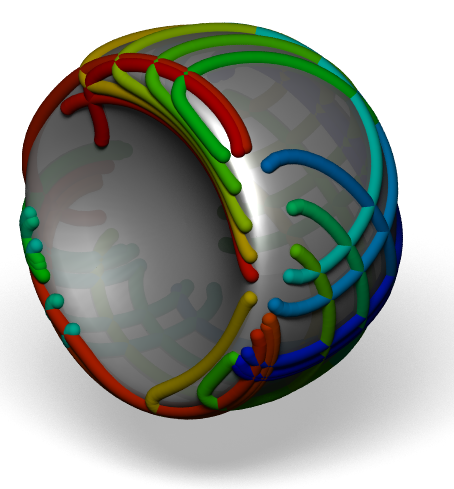}
    \put(27,55){\includegraphics[width=0.26\textwidth]{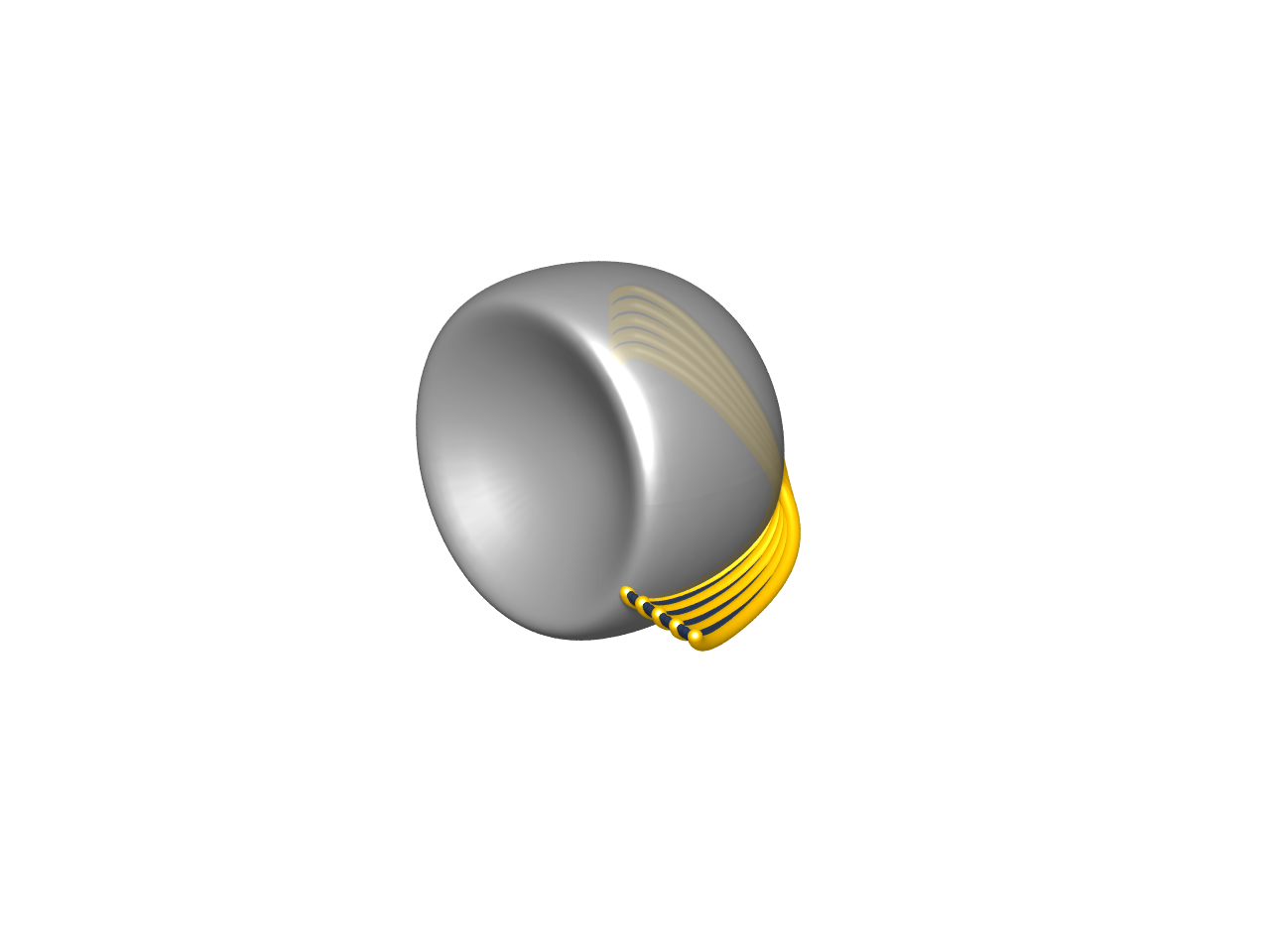}}
    \put(15,-30){\includegraphics[width=0.32\textwidth]{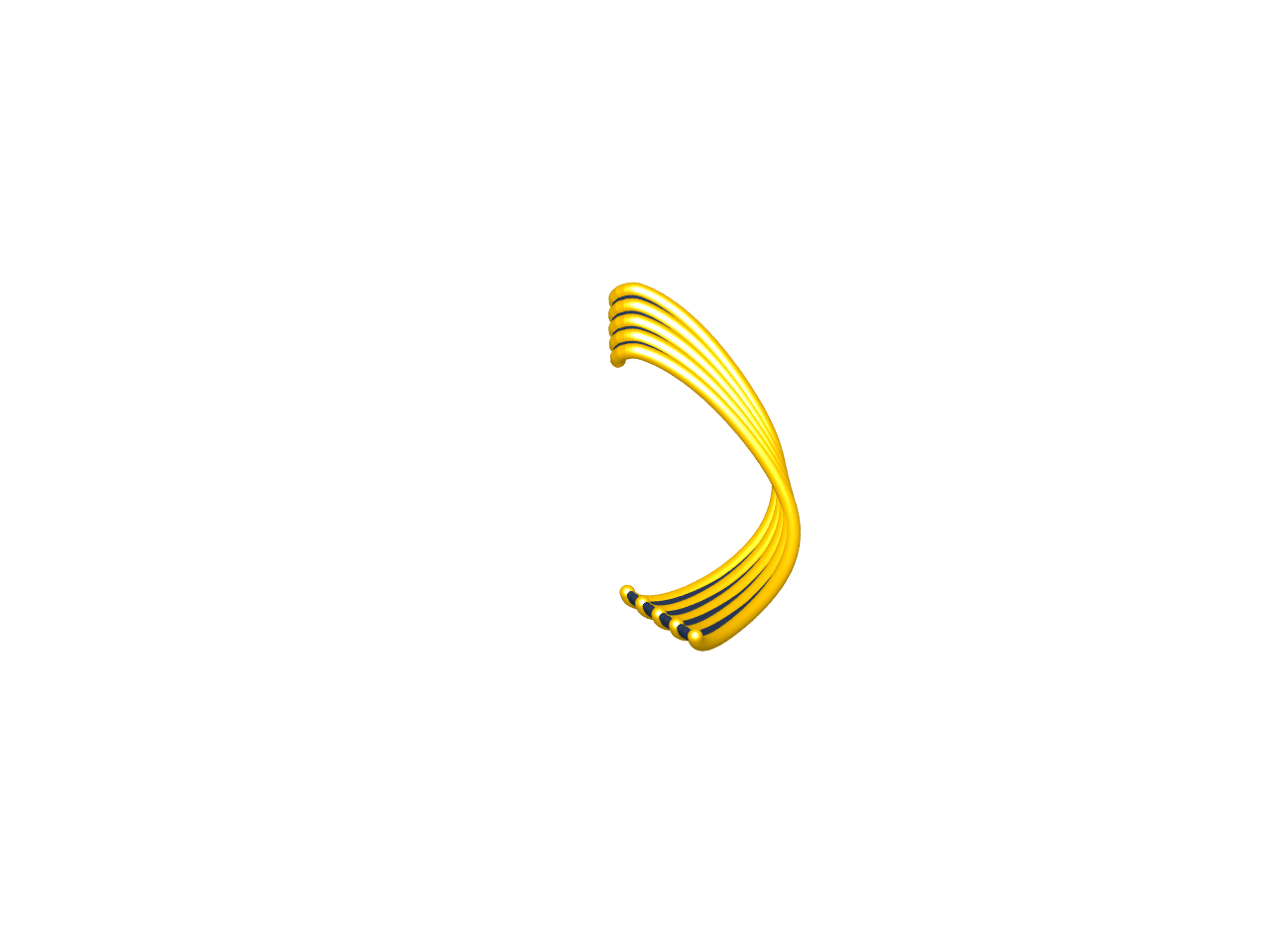}}
    \put(0,-5){(b) \quad $\bw=(0,0,0,1)$}
\end{overpic}}\hfill
\begin{overpic}[width=.2\textwidth]
    {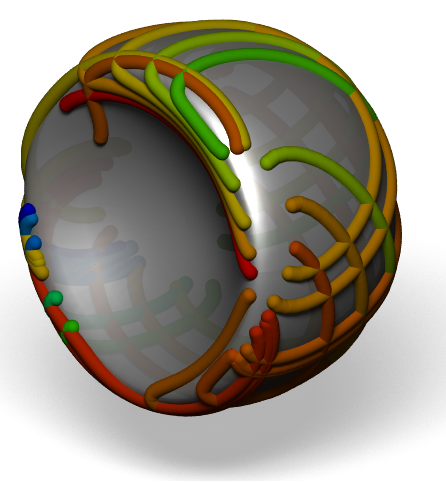}
    \put(30,55){\includegraphics[width=0.26\textwidth]{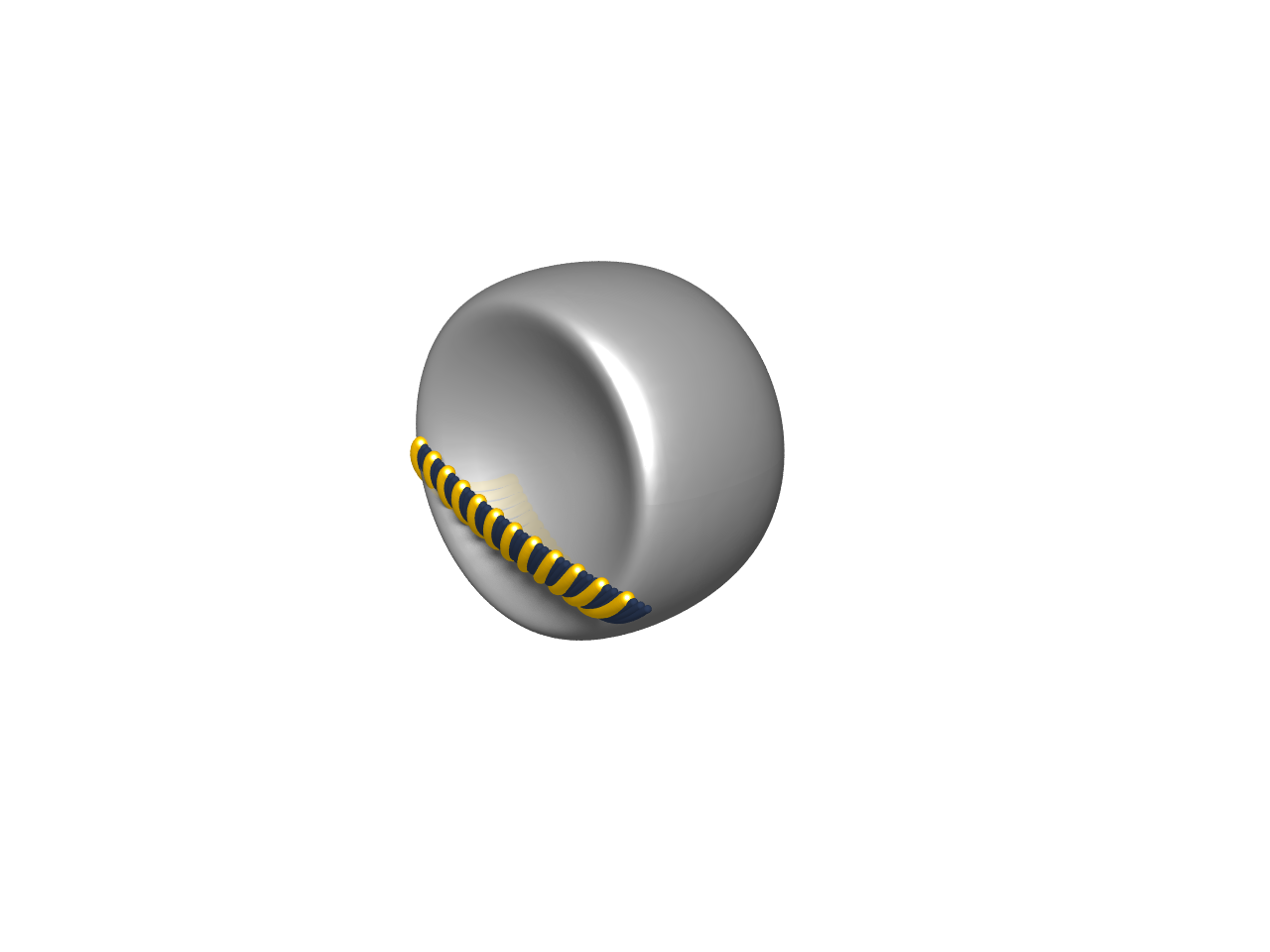}}
    \put(35,-30){\includegraphics[width=0.32\textwidth]{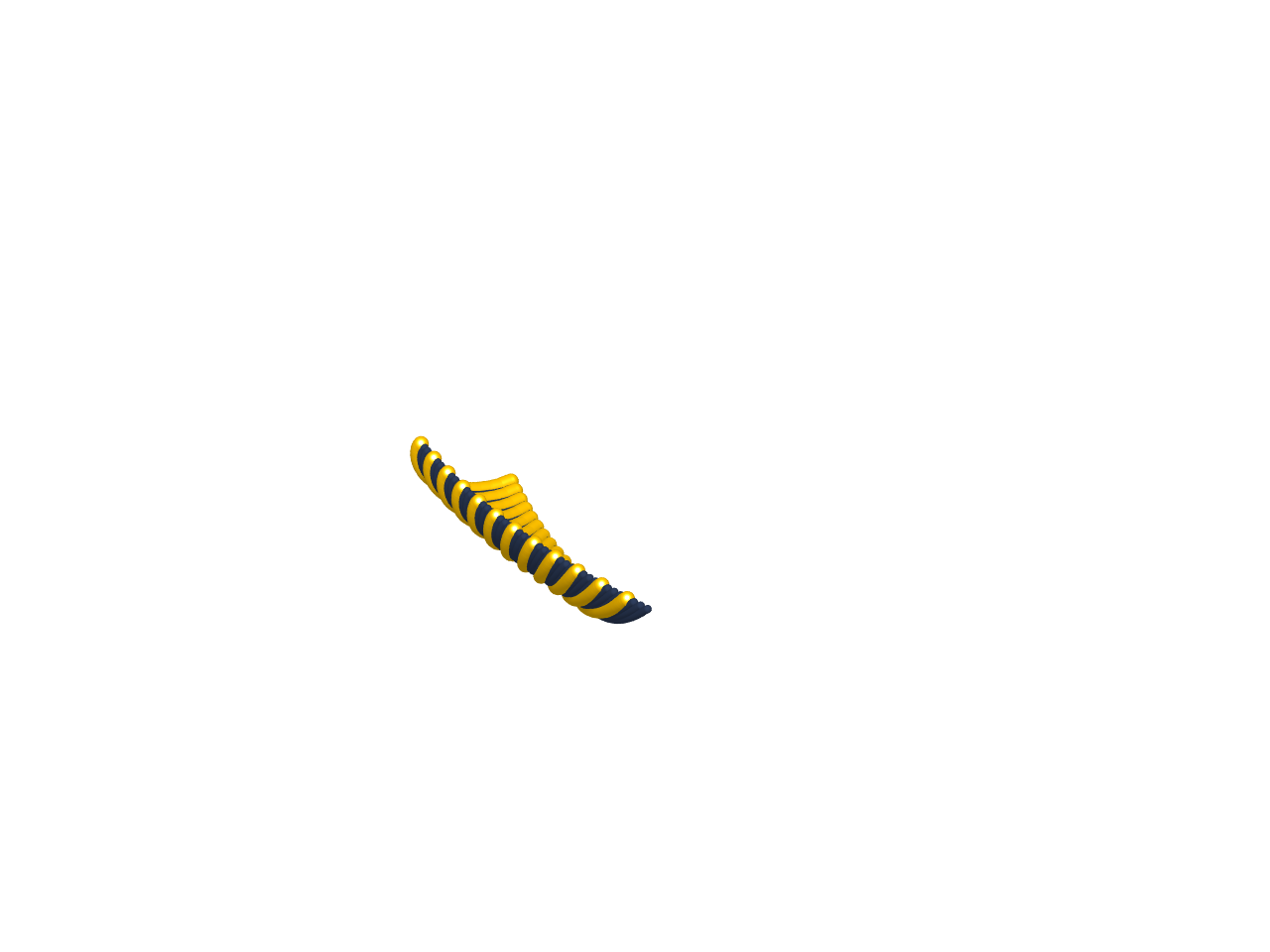}}
    \put(0,-5){(c) \quad $\bw=(0,1,0,0)$}
\end{overpic}\hfill
\raisebox{0.0\columnwidth}{
\begin{overpic}[width=.2\textwidth]
    {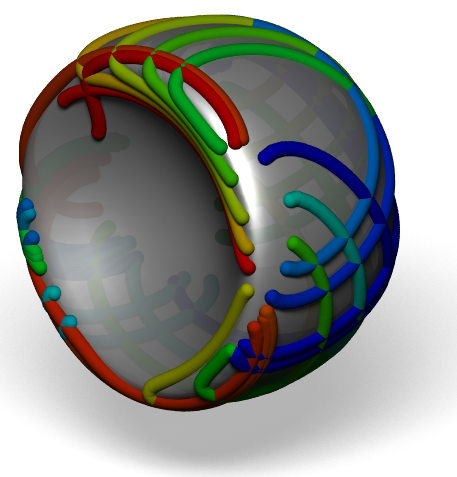}
    \put(28,55){\includegraphics[width=0.26\textwidth]{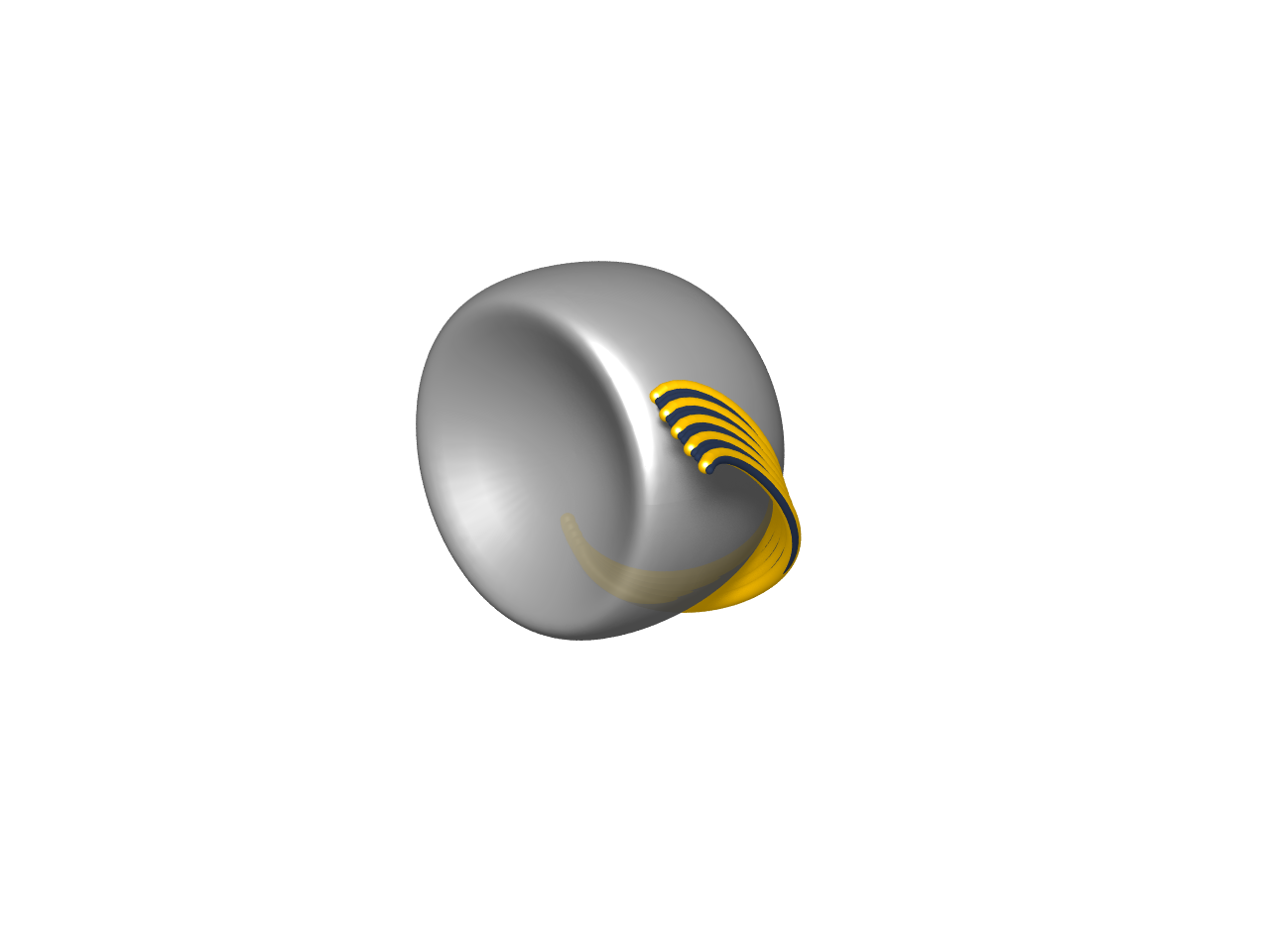}}
    \put(0,-40){\includegraphics[width=0.36\textwidth]{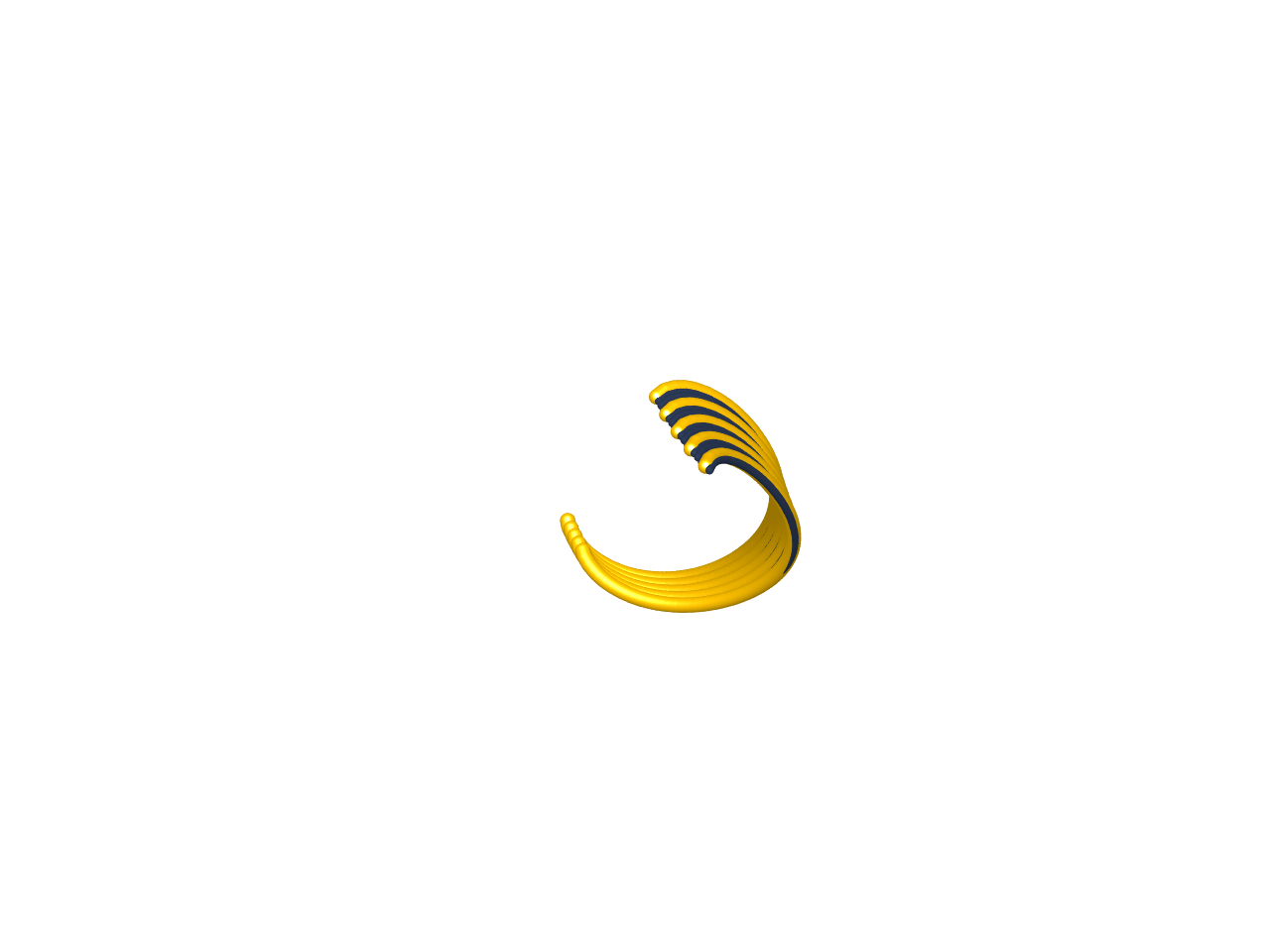}}
    \put(0,-5){(d) \quad $\bw=(0,0,1,0)$}
\end{overpic}}\hfill \vrule width0pt\\[-2ex]
\end{PDF}
    \caption{A family of first-order boundary curves from $\Gamma_1^\partial$, see (\ref{eq:family1}), on the smooth surface $\partial\Omega$.
    These curves are subsequently ranked by $E$ given in~\eqref{eq:E} considering several objectives: (b) the curves are additionally required to move as rigidly as possible, (c) as orthogonally as possible, (d) as parallel as possible. This is achieved by setting the appropriate weight vector $\bw$ in \eqref{eq:E}. The best seed curve and the corresponding stream surface are shown in each situation.}\label{fig:CrvdDomain}
\end{figure*}

\paragraph*{Alternative ranking energies.}
In the case when the strain energy is not the main objective, one may consider alternative components of a general \emph{ranking energy} $E$:
\begin{equation}\label{eq:E}
E(\gamma) = w_1 E_{in} + w_2 E_{ortho} + w_3 E_{para} + w_4 E_{rigid} + E_2,
\end{equation}
where
\begin{equation}\label{eq:AltEs}
\renewcommand{\arraystretch}{1.6}
\begin{array}{rcl}
E_{in}(\gamma) & = & \frac{1}{s_1-s_0}\int_{s_0}^{s_1}
1 - \langle \bs_t, \bm \rangle^2 \, \mathrm{d}s |_{t=0},\\
E_{ortho}(\gamma) & = & \frac{1}{s_1-s_0}\int_{s_0}^{s_1}
\langle \bs_t, \bs_s \rangle^2 \, \mathrm{d}s |_{t=0},\\
E_{para}(\gamma) & = & \frac{1}{s_1-s_0}\int_{s_0}^{s_1}
\langle \bs_{st}, \bs_{st} \rangle^2 \, \mathrm{d}s |_{t=0},\\
E_{rigid}(\gamma) & = & \frac{1}{s_1-s_0}
\int_{s_0}^{s_1} \langle \bar\bc + \bc \times \gamma(s), \bs_t \rangle^2 \, \mathrm{d}s |_{t=0}.\\
\end{array}
\end{equation}
The particular components of $E$ have the following meaning:
$E_{in}$ votes for curves that start moving to the interior part of $\Omega$, $\bm$ being the inward normal of $\partial\Omega$;
$E_{ortho}$ selects curves that start moving orthogonally to the flow, i.e., the tangent vectors $\bs_s=\dot{\gamma}$ are (in the least square sense) as orthogonal as possible to the velocity vectors $\bs_t=\bv$;
similarly $E_{para}$ finds curves that propagate in a parallel fashion; and
$E_{rigid}$ favours curves that move as rigidly as possible,
$(\bc, \bar\bc)$ being the instantaneous motion; see
\cite{Pottmann-2001-clg,Barton-2012-Snakes}.

The behaviour of the ranking energy $E$ depending on weights $\bw = (w_1, w_2, w_3, w_4)$ when applied on a curved domain $\Omega$ is shown in Fig.~\ref{fig:CrvdDomain}.

However, there exist rare scenarios with no first-order boundary curves, i.e., the set $\Gamma_1^\partial$ is empty, or the existing ones are insufficient for a particular application. If that is the case, we turn to our third strategy, which is always guaranteed to produce first-order curves.

\subsection{First-order interior curves}\label{sec:str}
In the rare situation when there are no second-order curves ($\Gamma_2 = \emptyset$) and no first-order boundary curves ($\Gamma_1^\partial = \emptyset$), we identify first-order curves that are, given an initial point and first-order direction, curvature minimising among those in $\Gamma_1$. The benefit here is twofold: such curves are always guaranteed to exist, and they are particularly well suited for visualisation since they are as straight as possible.
These curves are given by point-wise minimisation of $||\ddot{\gamma}(s)||$ subject to $\dot{\gamma}(s)\bj(\gamma(s))\dot{\gamma}(s)^\top=0$; cf.~\eqref{eq:conjugate}.

Discretisation and implementation details are discussed in Section~\ref{sec:impl}. Before all that, we proceed to the exploration of optimal stream surface generation.

\section{Stream surfaces and optimisation}\label{ssec:StrainMinSrf}

We start by defining strain minimising stream surfaces, which are then optimised with respect to a certain strain energy.

\subsection{Strain minimising stream surfaces}\label{sec:smss}

Second-order curves, or first-order curves if the former ones do not exist, are good initial guesses for seed curves, called \emph{candidate seed curves},
but only \emph{locally}. Our goal is to find stream surfaces for which their generating seed curves propagate \emph{globally} in a certain strain-minimising manner. This is formalised in
\begin{dfn}
The \emph{strain energy} of a stream surface $\bs$ is given by
\begin{equation}\label{eq:Esrf}
E_{\bs} = \frac{1}{\textnormal{area}(\bs)}\int_{t_0}^{t_1} E_1(\bs(s,t)) \, \mathrm{d}t.
\end{equation}
A stream surface that minimises this energy will be called \emph{strain minimising}.
\end{dfn}
In other words, $E_{\bs}$ measures the strain given by the deformation of the seed curve of $\bs$ through the field by
accumulating the strain energies of all its timelines.
This energy is, as in the case with seed curves, used as a ranking criterion for stream surfaces.

\subsection{Optimisation of stream surfaces}\label{sec:opt}

Candidate seed curves generate stream surfaces with relatively low strain energy. However, their effect is still only \emph{local}, whilst
we are interested in minimising the strain energy \emph{globally}.
This is achieved by employing an optimisation procedure; see Fig.~\ref{fig:flowchart} for a schematic overview.

\begin{figure*}[!tbh]
\begin{PDF}
 \hfill
 \begin{overpic}[width=.19\textwidth]
	{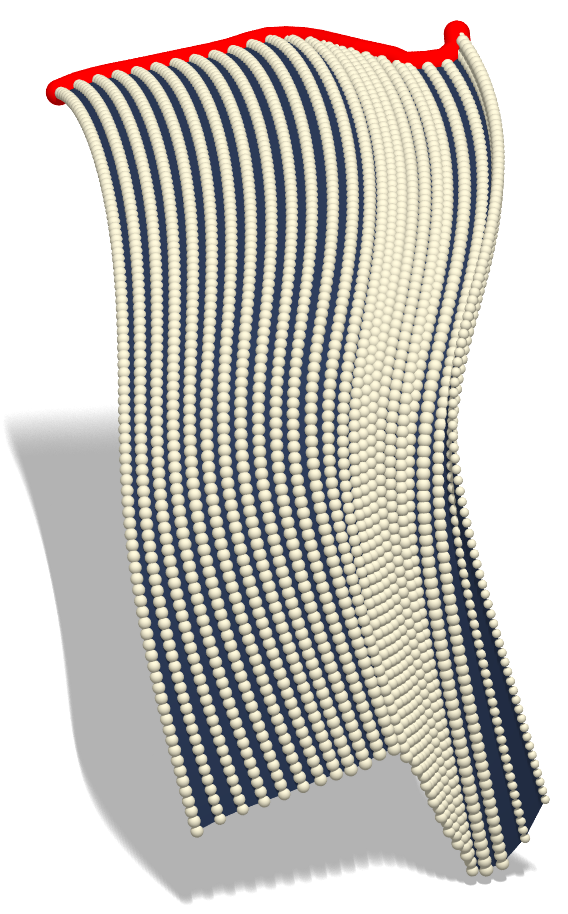}
	\put(2,95){$\gamma_0$}
    \put(15,1){\tiny{$E_{\bs_0} = 7.41 \cdot 10^{-4}$}}
    \put(0,0){(a)}
	\end{overpic}\hfill
  \begin{overpic}[width=.19\textwidth]
	{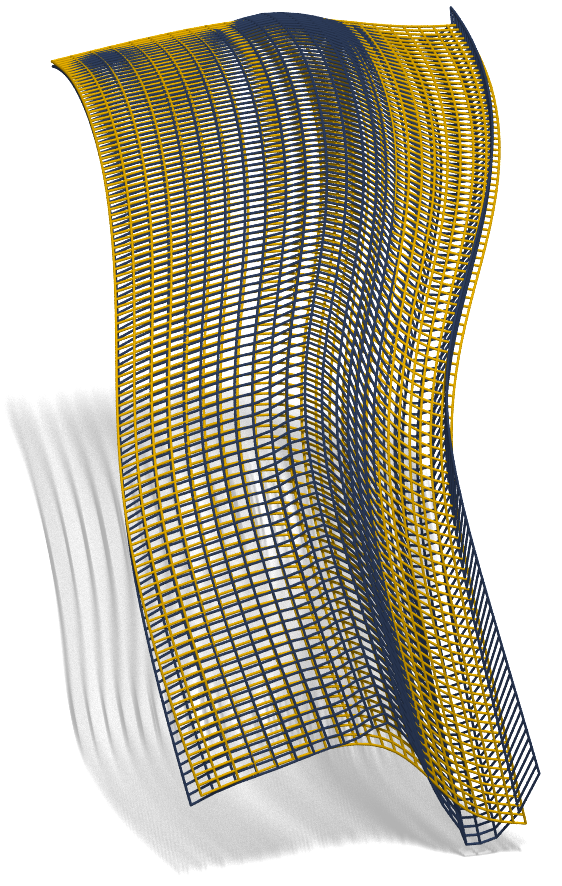}
	\put(0,0){(b)}
	\end{overpic}\hfill
 \begin{overpic}[width=.19\textwidth]
	{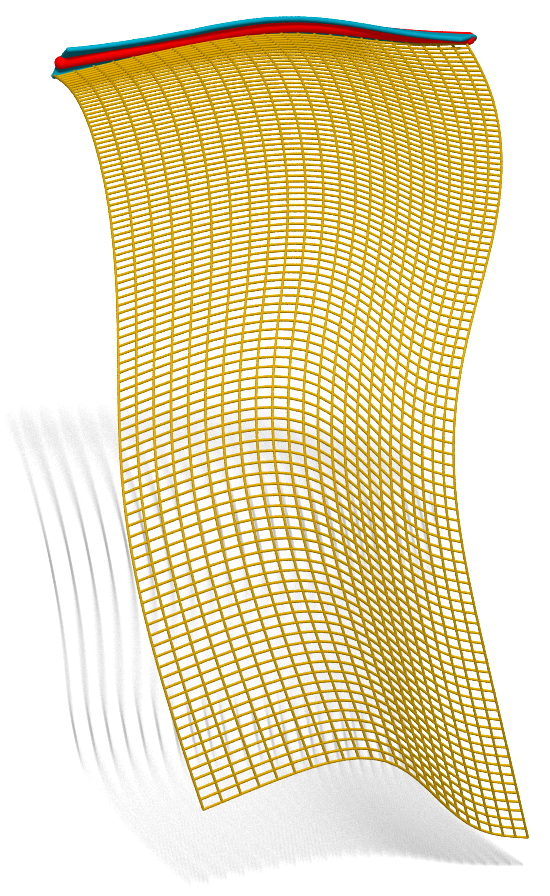}
    \put(0,0){(c)}
    \put(0,95){\color{cyan}$\mathcal S_1$}
    \put(40,2){$\bs_0^{opt}$}
	\end{overpic}\hfill
\begin{overpic}[width=.19\textwidth]
	{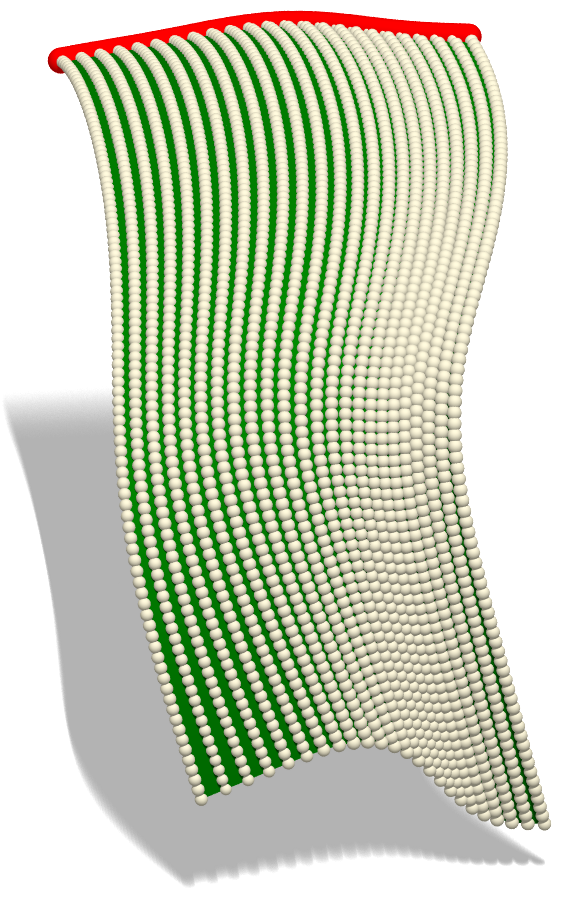}
    \put(2,98){$\gamma_1$}
	\put(15,1){\tiny{$E_{\bs_1} = 2.86 \cdot 10^{-4}$}}
    \put(0,0){(d)}
	\end{overpic}\hfill
\begin{overpic}[width=.19\textwidth]
	{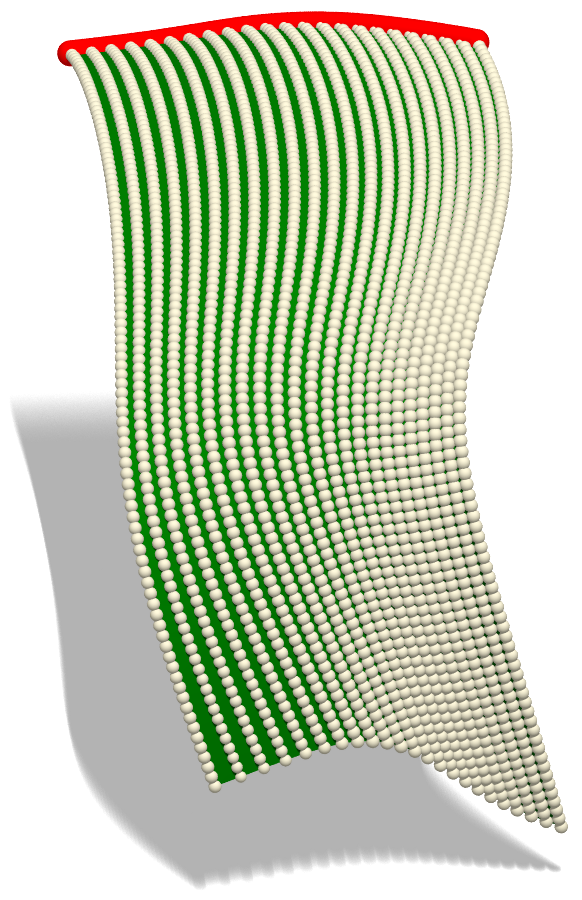}
    \put(2,98){$\gamma_2$}
	\put(15,1){\tiny{$E_{\bs_2} = 2.32 \cdot 10^{-4}$}}
    \put(0,0){(e)}
	\end{overpic}\hfill
\vspace{-2pt}
\end{PDF}
  \Acaption{1.5em}{Stream surface optimisation.
	(a) The seed curve $\gamma_0$ (red) is integrated, resulting in the initial stream surface $\bs_0$ (blue). The radii of the spheres correspond to arc-length change.
	(b) The stream surface is discretised into a quad mesh $\bQ$ and optimised.
	(c) Since the optimised surface $\bs_0^{opt}$ is not a stream surface any more, its timelines are back-integrated to $t=0$, forming a set $\mathcal S_1$ from which the original seed curve is improved to $\gamma_1$ (red).
	(d) The improved seed curve gives rise to an updated stream surface $\bs_1$. (e) The stream surface $\bs_1$ after another iteration of the optimisation cycle.
  }\label{fig:Opt}
  \end{figure*}

Candidate curves are used for initialising an optimisation cycle, which works as follows; see Fig.~\ref{fig:Opt}. A candidate seed curve $\gamma_0$ at $t=0$ is integrated to form its stream surface $\bs_0$,
which is subsequently optimised with respect to its strain energy \eqref{eq:Esrf}. The optimised surface $\bs_0^{opt}$, however, is not, in general, a stream surface any more. Therefore, its timelines are back-integrated to the initial time $t=0$ to form a set of space curves $\mathcal S_1$.
These curves are used to compute an updated $\gamma_1$ by least square fitting, which is forward integrated to build $\bs_1$ and so on.
The algorithm continues updating $\bs_i \gets \bs_{i+1}$ until the surface strain energy $E_{\bs_{i}}$ stops being improved,
or when the maximum number of iterations, set to ten if not stated otherwise, is reached.
The particular steps of the algorithm are explained, including implementation details, in the next section.
%

\section{Implementation}\label{sec:impl}

\begin{figure*}[t]
\begin{center}
\fcolorbox{olive}{white}{\begin{minipage}{.87\textwidth}
{\footnotesize\color{blau}
 \tikzstyle{block} = [rectangle, draw, fill=drot!5,
    text width=6.3em, text centered, rounded corners, minimum height=2em]
 \tikzstyle{info} = [rectangle, color=black, text width=12em, minimum height=2em]
 \tikzstyle{line} = [draw, -latex']
\begin{tikzpicture}[node distance=8.2em]
    \node [block] at (0,0) (cuts)
	{Compute first- or second-order curves, Sections~\ref{sec:SE1} and~\ref{sec:SE2}};
    \node [block, right of=cuts] (bins)
	{Integrate the vector field, obtain initial stream surfaces};
    \node [block, right of=bins] (ribs)
	{Rank the seed curves according to $E_\bs$, Eq.~(\ref{eq:Esrf})};
    \node [block, right of=ribs] (hermite)
	{For each seed curve, optimise its stream surface, Section~\ref{sec:opt}};
    \node [block, right of=hermite] (profile)
	{Back-integrate, improve seed curve, see Fig.~\ref{fig:Opt}(c)};
    \node [block, right of=profile] (measure)
	{Check if stream surface is strain minimising};
    \node [block, right of=measure] (done)
	{Return stream surface};
  \node [info] at (9.6,-1.0) {If $E_{\bs}$ is not good enough};
  \node [info] at (2.8,-1.0) {Go to regions with low $E_\bs$};
    \path [line] (cuts) -- (bins);
    \path [line] (bins) -- (ribs);
    \path [line] (ribs) -- (hermite);
    \path [line] (hermite) -- (profile);
    \path [line] (profile) -- (measure);
    \path [line] (measure) -- (done);
    \draw[line] (measure) .. controls  (12,-1.5) and (7,-1.5) .. (hermite);
    \draw[line] (ribs) .. controls  (5,-1.5) and (0,-1.5) .. (cuts);
 \end{tikzpicture}}
 \vspace*{-1.7em}
\begin{minipage}[b]{.99\textwidth}\caption{\label{fig:flowchart} Algorithm overview.}
\end{minipage}
\end{minipage}}
\end{center}
\end{figure*}
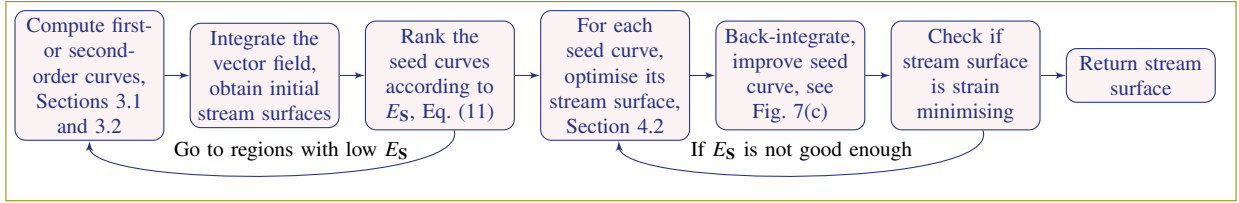

Having explored the smooth setting of seed curves and stream surfaces, we now discretise these notions and discuss the implementation
of our algorithm; see Fig.~\ref{fig:flowchart} for an overview.

\paragraph*{Finding candidate curves.} Our approach is based on (adaptive) sampling. In the case of first-order boundary curves, $\partial\Omega$ is explored; see examples in Figs.~\ref{fig:EnergyRanking} and \ref{fig:CrvdDomain}. Otherwise, the interior of $\Omega$ is sampled and first- and second-order curves are computed. These are ranked by \eqref{eq:Esrf} and the regions with low values of $E_\bs$ are sampled with higher density. This sampling is repeated recursively, if not stated otherwise, three times and the best $5\%$ are taken as candidate curves.

Seed curve computation proceeds as follows. Starting with a sampled seed point, the admissible tangent field is computed (see Sections \ref{sec:SE2}, \ref{sec:bnd}, and \ref{sec:str} for particular cases) and integrated in all admissible directions. The step-size is, by default, set as $1\%$ of the domain's diameter $\mathrm{diam}(\Omega)$, but can be refined if the sampling is not sufficient; see Fig~\ref{fig:Turbine}.
Only curves of length greater or equal to $\mathrm{diam}(\Omega)/10$ are considered as candidates. The integration terminates if there is no admissible direction to continue in, the maximum length set to $\mathrm{diam}(\Omega$) is achieved, or, for non-boundary curves, the boundary is reached. Integration of first-order interior curves can be found in Appendix C.

\paragraph*{Integrating the vector field.}
Various techniques for numerical \emph{streamline} integration
have been studied \cite{Nielson-1997-SV}. Higher-order methods
with an adaptive stepsize were shown to be the preferable choice when considering
accuracy versus speed trade-offs. In our case, however,
since the back-integrated timelines of the optimised stream surface are taken as input for seed curve optimisation, our prime focus is on accuracy. Therefore,
we employed the classical fourth order Runge-Kutta method with constant stepsize.

\paragraph*{Optimising the stream surface.}
Given a stream surface $\bs$, we want to reduce its strain energy \eqref{eq:Esrf}.
Sampling $m$ points in the $s$-direction (seed curve) and $n$
points in the $t$-direction (time), $\bs$ is discretised into a quad mesh $\bQ$ having $m \times n$ vertices $\bq_{i,j}^0$.
Denoting the vertices of the optimised surface by $\bq_{i,j}$ and setting $\be_{i,j} := \bq_{i+1,j+1} - \bq_{i+1,j} + \bq_{i,j} - \bq_{i,j+1}$,
the desired improvement is formulated as a non-linear least squares optimisation with the objective function
\begin{equation}\label{eq:F}
F(\bQ) = F_{strain}(\bQ) + \mu_1 F_{fair}(\bQ) + \mu_2 F_{prox}(\bQ)
\, \, \textnormal{with}
\end{equation}
\vspace*{-15pt}
\begin{equation}\label{eq:Fcomp}
\renewcommand{\arraystretch}{1.6}
\begin{array}{rcl}
F_{strain}(\bQ) & = & \sum_{i,j}\langle \bq_{i+1,j} - \bq_{i,j}, \be_{i,j}\rangle^2\\
F_{fair}(\bQ) & = & \sum_{i,j} \| \bq_{i+1,j} -2\bq_{i,j} + \bq_{i-1,j}\|^2+\\
 & & + \| \bq_{i,j+1} -2\bq_{i,j} + \bq_{i,j-1}\|^2,\\
F_{prox}(\bQ) & = & \sum_{i,j} \| \bq_{i,j} - \bq_{i,j}^0\|^2,
\end{array}
\end{equation}
where $F_{strain}$ reflects the strain minimising condition (cf.~\eqref{eq:Esrf}), $F_{fair}$ is a fairness term, and $F_{prox}$ is a proximity term,
a regulariser that forces the vertices of the optimised mesh not to deviate much from the input. The optimisation problem is solved using the Gauss-Newton method for all the examples in the paper and the accompanying video. Experimentally, the weights were set to $\mu_1 = 0.1$, $\mu_2 = 0.02$. This results in an optimised surface $\bs^{opt}$.

\paragraph*{Improving the seed curves.}
Having obtained the optimised surface $\bs^{opt}$, its timelines $\bs^{opt}(s,t_j) = \gamma_j$, $j=1,\dots,n$ are back-integrated to the initial time
instant $t=0$, giving the set $\mathcal S_1$ of space curves $\gamma_j^{t=0}$, see Fig.~\ref{fig:Opt}(c). If $\bs^{opt}$ was an exact stream surface, all $\gamma_j^{t=0}$ would coincide.
We improve the seed curve $\gamma$ of $\bs$ by replacing it by the least squares approximation \cite{Hoschek-2002-CAGD} of $\gamma_j^{t=0}$, $j=1,\dots,n$.
Note that we have the information about the correspondence in the $s$-direction, i.e., for a fixed $i$, all the points $\bq_{ij}$, $j=1,\dots n$, need to correspond to a single point $\gamma(s_i)$. This fact simplifies the problem to point-wise averaging.
Having numerically back-integrated $\bs^{opt}$ in a point-wise fashion resulting in $\bq_{ij}^{t=0}$, the seed curve update is achieved by setting $\gamma(s_i) = \frac{1}{n}\sum_j\bq_{ij}^{t=0}$.

\paragraph*{Curve and surface trimming.} So far, we have not discussed how to set the intervals $[s_0,s_1]$ and $[t_0,t_1]$. In practise, both of these have to be finite.
In the $s$ (seed curve) direction, the tangent vector field is integrated while admissible directions exist, or until the boundary is reached.
In the $t$ (time) direction, the integration is terminated if the timeline reaches the boundary, or, as shown in Fig.~\ref{fig:Rayleigh}, when the maximum number of timesteps is reached.

\section{Numerical results}\label{sec:exmp}

\begin{figure}[!tbh]
\begin{PDF}
 \hfill
 \begin{overpic}[width=.95\columnwidth]
	{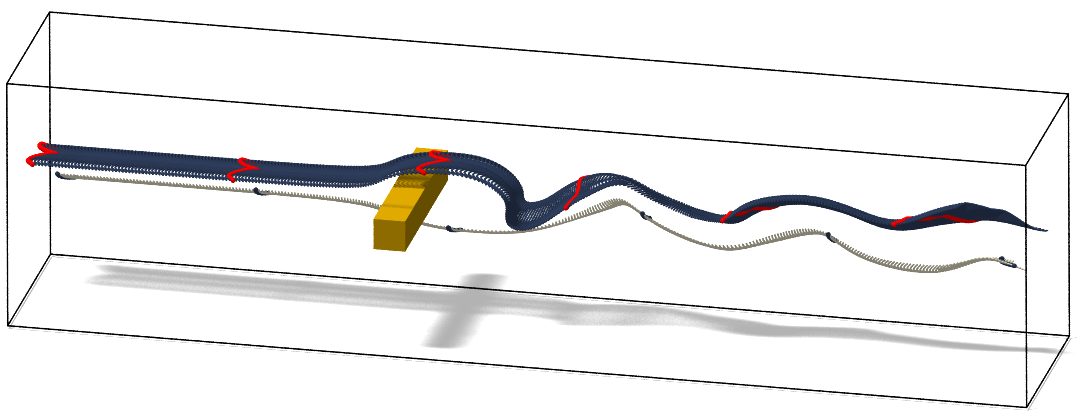}
    \put(0,0){(a)}
	\end{overpic}\hfill
\vspace{-0.01cm}
 \begin{overpic}[width=.49\columnwidth]
	{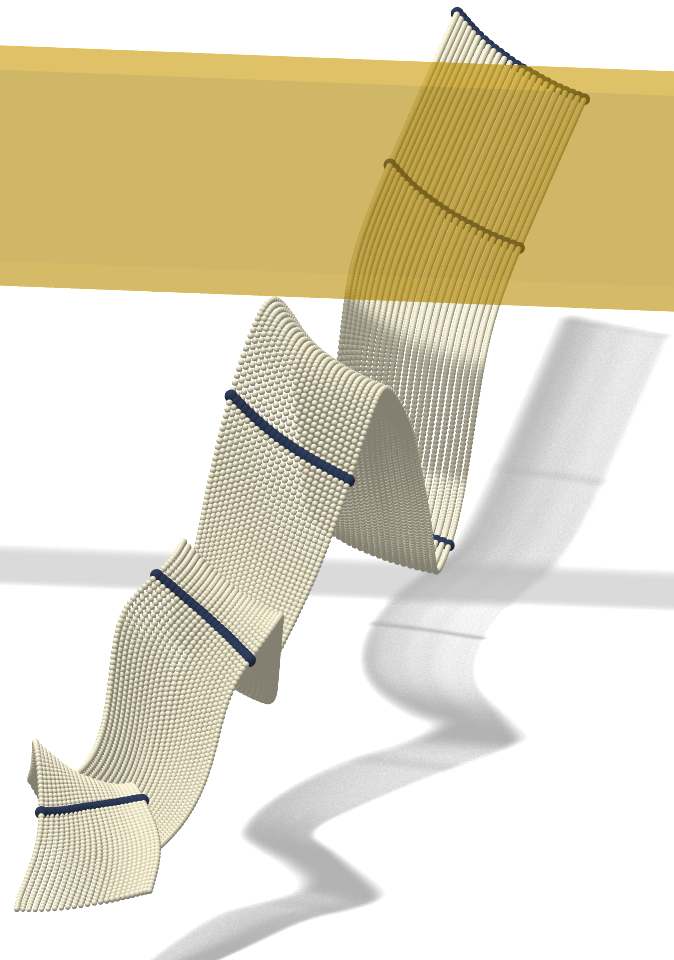}
     \put(50,0){(b)}
	\end{overpic}\hfill
\begin{overpic}[width=.49\columnwidth]
	{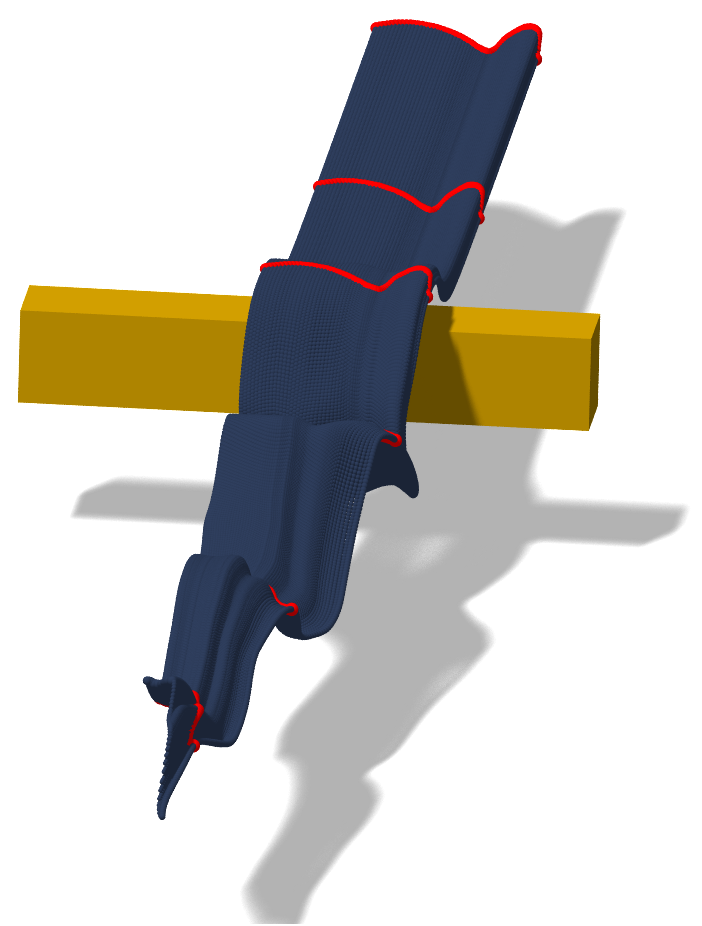}
	 \put(50,0){(c)}
	\end{overpic}\hfill
\vspace{-2pt}
\end{PDF}
  \Acaption{1.5em}{The `square cylinder' of \cite{camarri05}. (a) Incompressible flow moving from left to right in a box, encountering a rigid obstacle (yellow box). The best strain-minimising surfaces for seed curves of different lengths are shown. The length is demanded to be at least $5\%$ (b) and $10\%$ (c) of the domain's diameter. }\label{fig:SquareCylinder}
  \end{figure}

The algorithm was tested on several benchmark datasets. The example in Fig.~\ref{fig:SquareCylinder} is a direct numerical Navier-Stokes simulation by \cite{camarri05} that is publicly available \cite{iCFDDatabase}. We used a uniformly resampled version, which has been provided by Tino Weinkauf and used in \cite{Weinkauf-2008-SquareCylinder}. The example is based on the last time instant of the unsteady flow.

\begin{figure}[!tbh]
\begin{PDF}
    \hfill
        \begin{overpic}[width=.49\textwidth]
    	{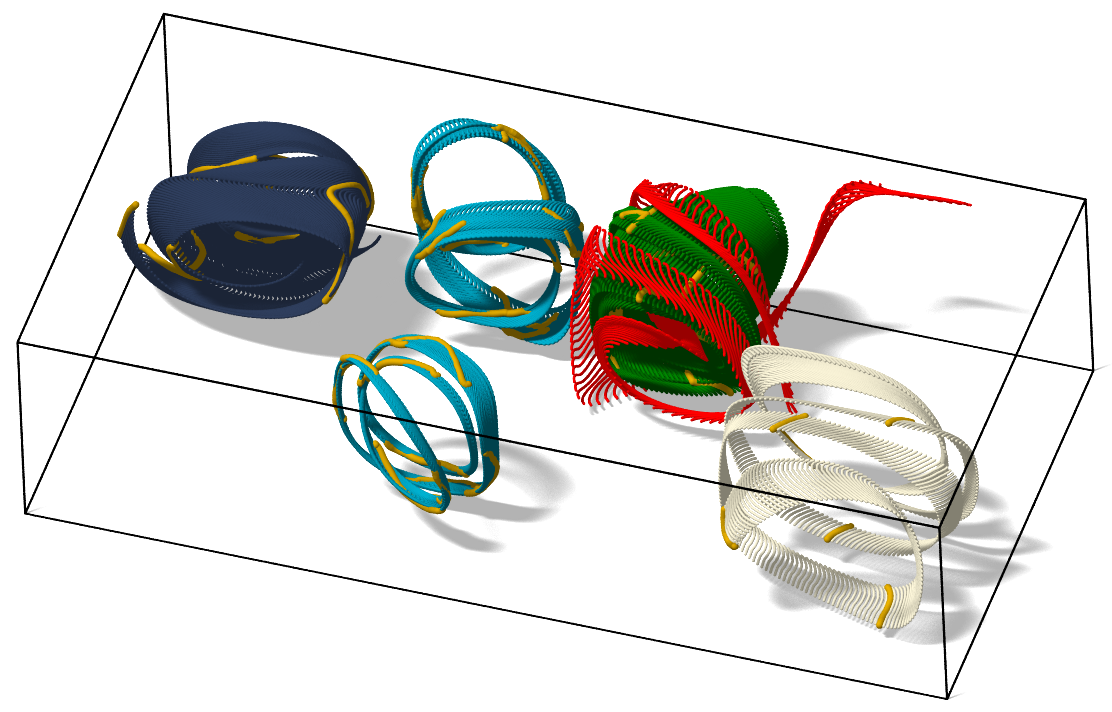}
        \put(75,51.5){\includegraphics[width=0.12\textwidth]{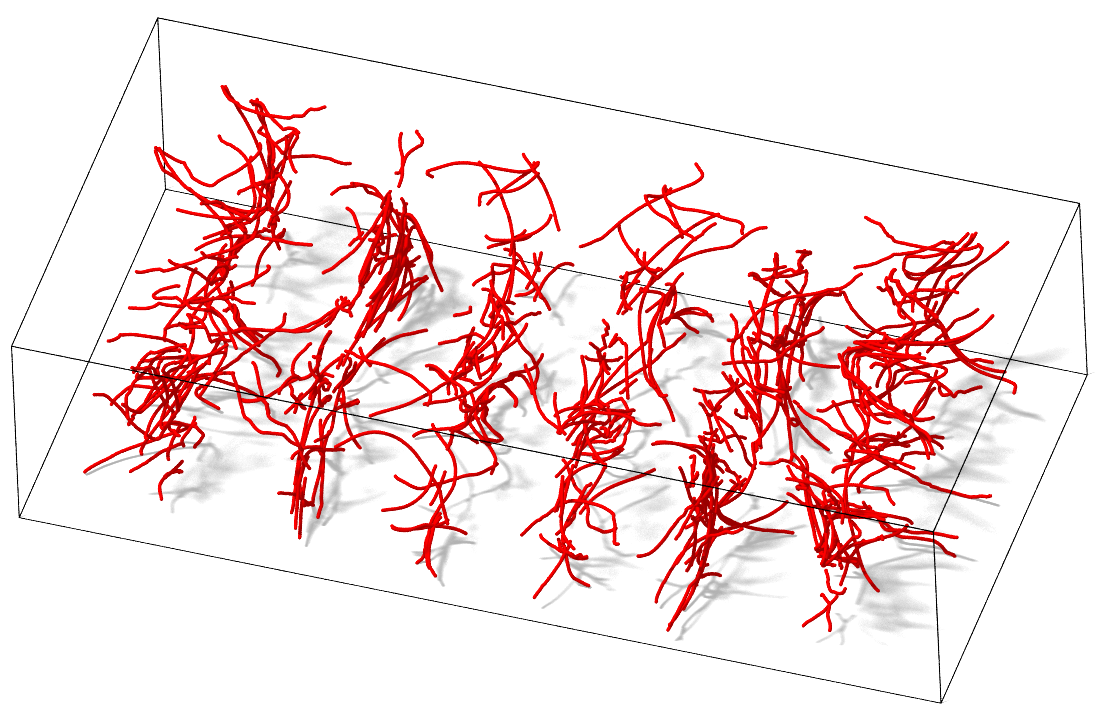}}
        \put(2,10){(a)}
        \put(0,-80){\includegraphics[width=.49\columnwidth]{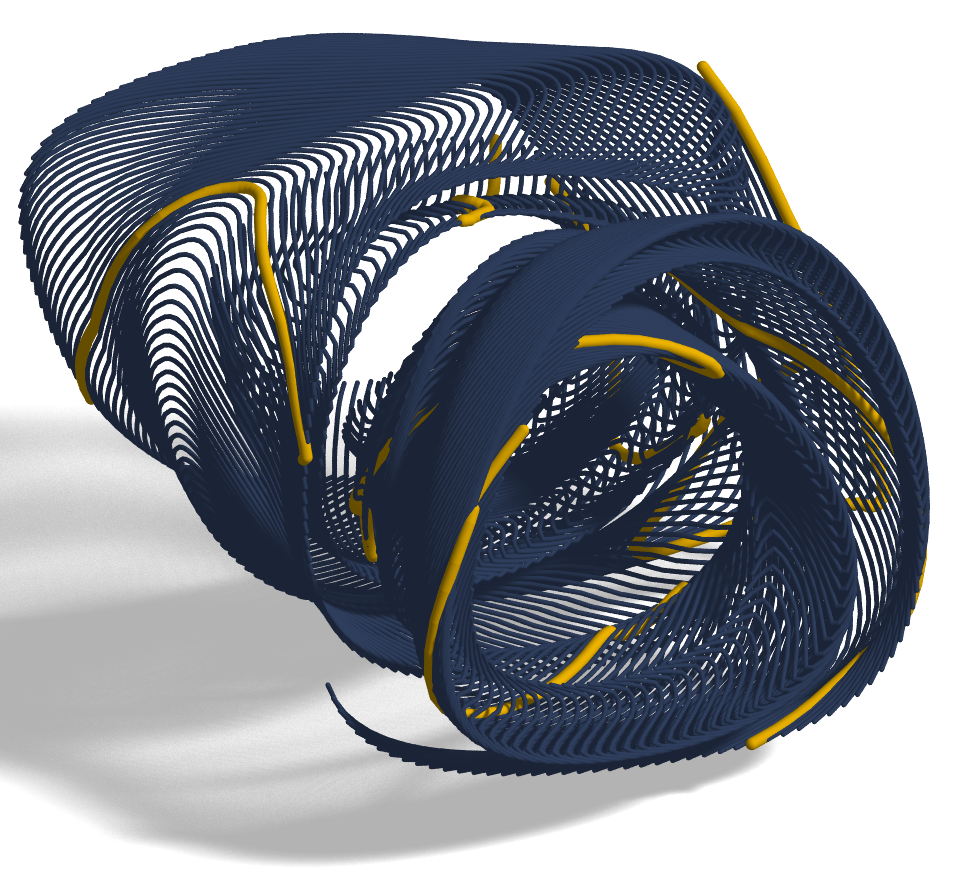}}
        \put(2,-80){(d)}
        \put(50,-80){\includegraphics[width=.43\columnwidth]{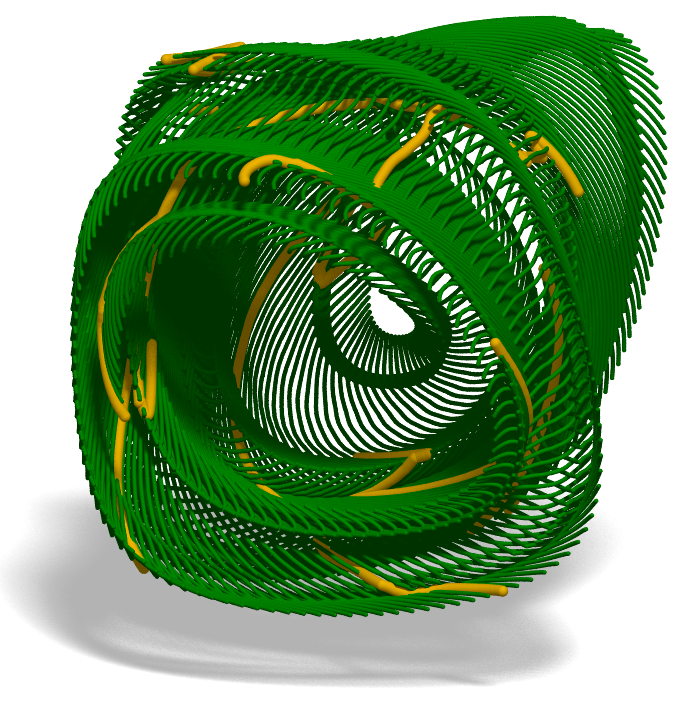}}
        \put(52,-80){(e)}
        \put(0,-35){\includegraphics[width=.58\columnwidth]{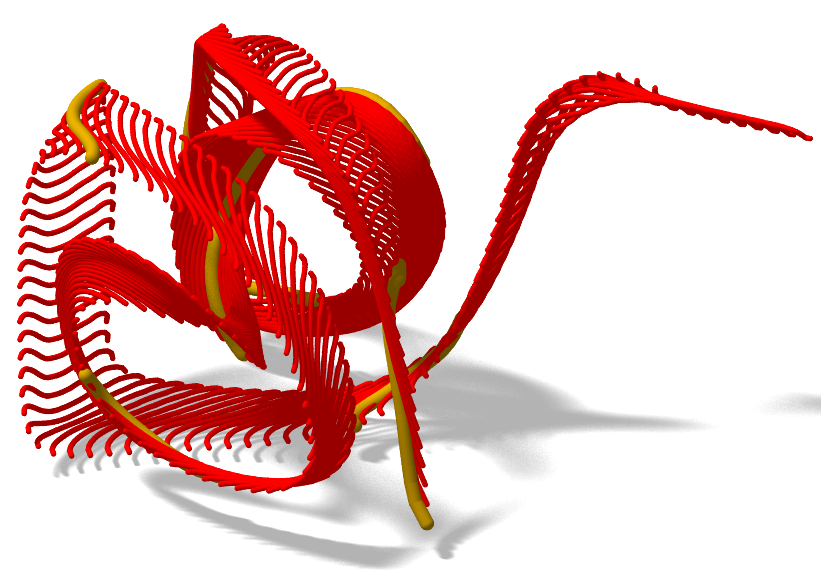}}
        \put(2,-35){(b)}
        \put(50,-35){\includegraphics[width=.39\columnwidth]{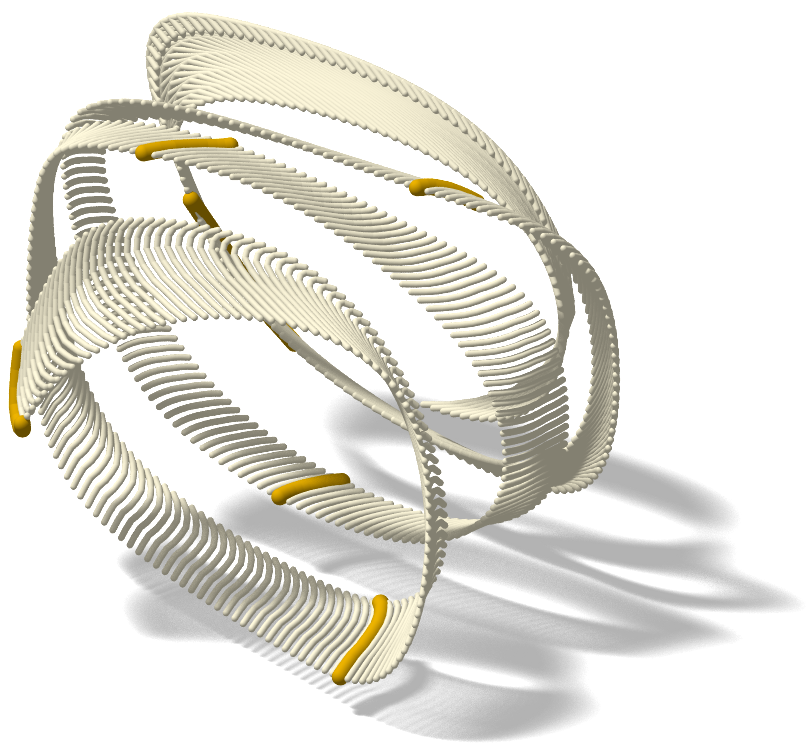}}
        \put(52,-35){(c)}
    	\end{overpic}
    \hfill
    \vspace{6.5cm}
\end{PDF}
  \Acaption{1.5em}{Rayleigh-B\'{e}nard heat convection. The dataset used here is given by one time instant of an unsteady vector flow obtained from \cite{NaSt3DGP}. (a) Several best strain minimising stream surfaces detected by our algorithm (see Fig.~\ref{fig:flowchart}) are shown. A total of $434$ second order curves were used for initialisation (top red); the optimisation parameters were set to $\mu_1 = 0.1$, $\mu_2 = 0.02$; see~\eqref{eq:F}.
  The vector field integration was terminated by (b--d) reaching the boundary and (e) by exceeding the upper bound on the number of timesteps.
}\label{fig:Rayleigh}
\end{figure}

\begin{figure*}[!tbh]
\begin{PDF}
  \hfill
 \begin{overpic}[width=.33\textwidth]
	{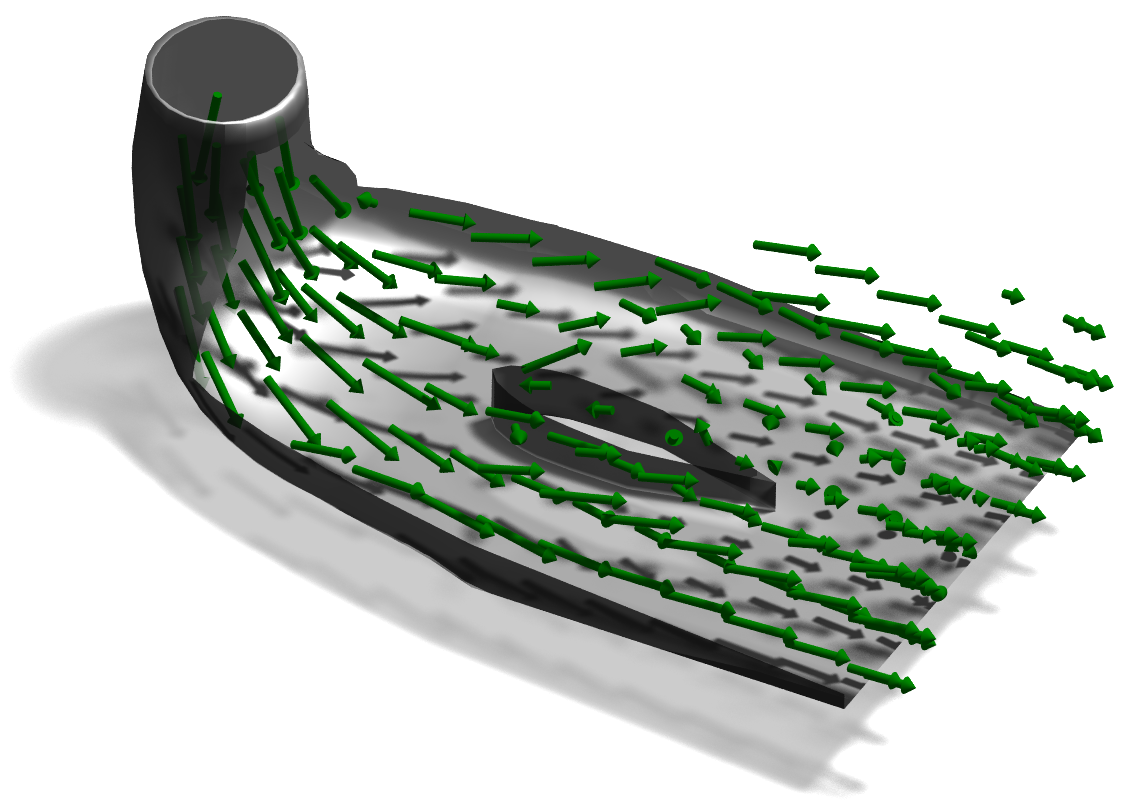}
    \put(0,5){(a)}
	\end{overpic}\hfill
  \begin{overpic}[width=.33\textwidth]
	{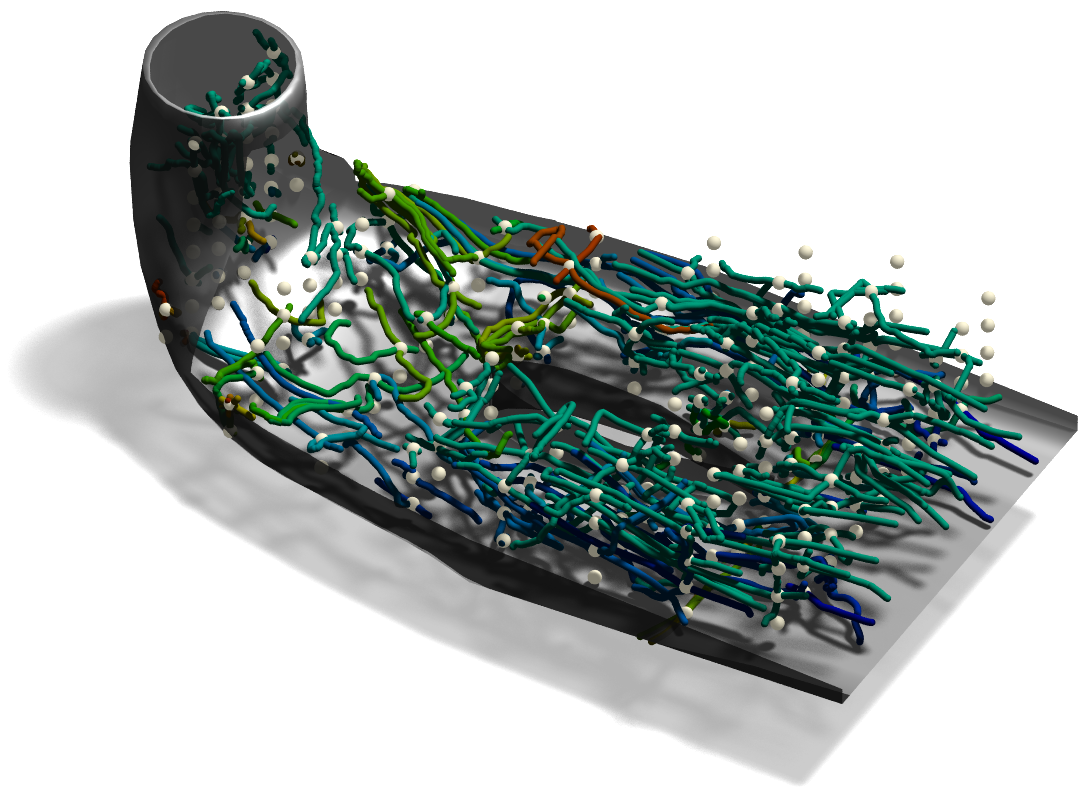}
	\put(0,5){(b)}
	\end{overpic}\hfill
 \begin{overpic}[width=0.33\textwidth]
	{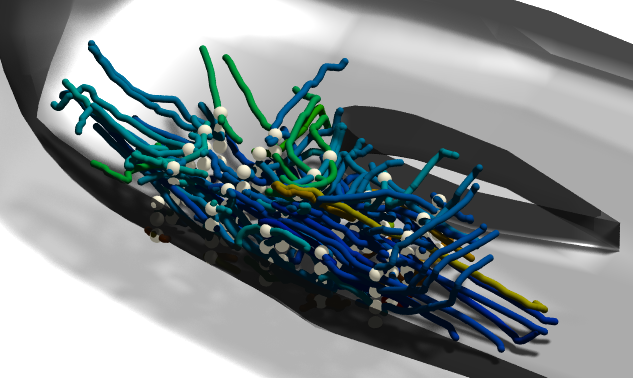}
	\put(0,5){(c)}
	\end{overpic}\hfill
\vspace{-0.01cm}
 \hfill
 \begin{overpic}[width=.24\textwidth]
	{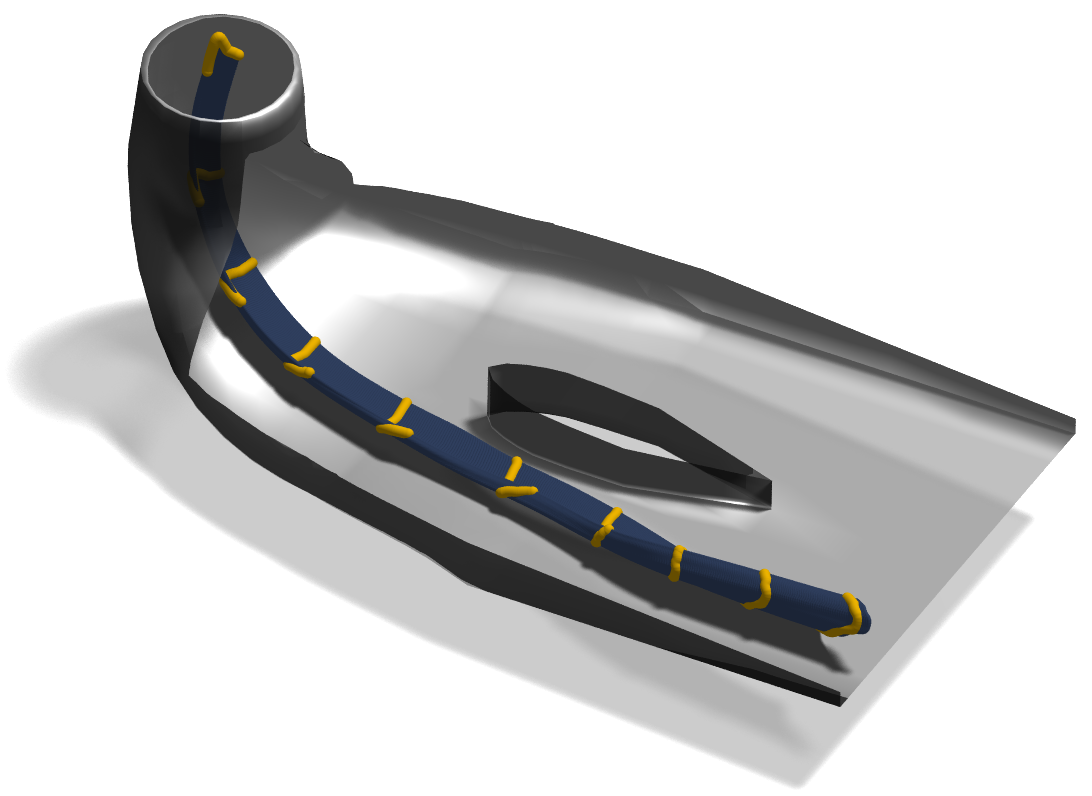}
	\put(0,5){(d)}
	\end{overpic}\hfill
  \begin{overpic}[width=.24\textwidth]
	{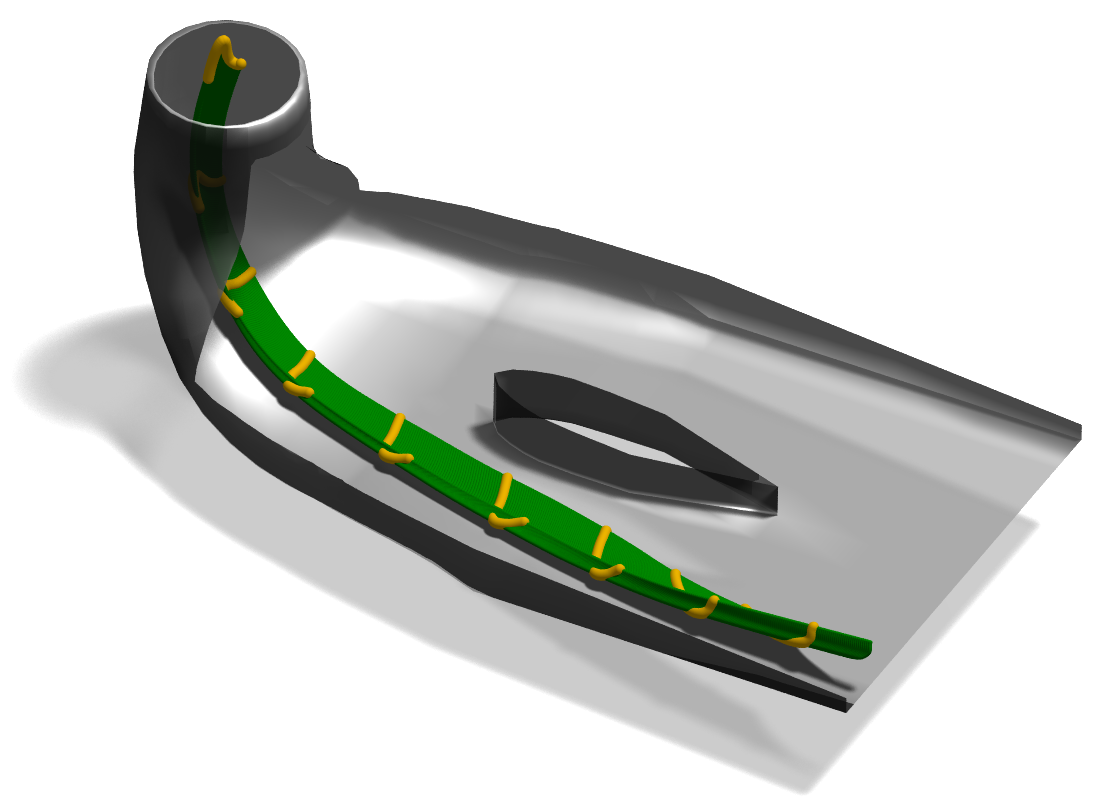}
	\put(0,5){(e)}
	\end{overpic}\hfill
 \begin{overpic}[width=.24\textwidth]
     {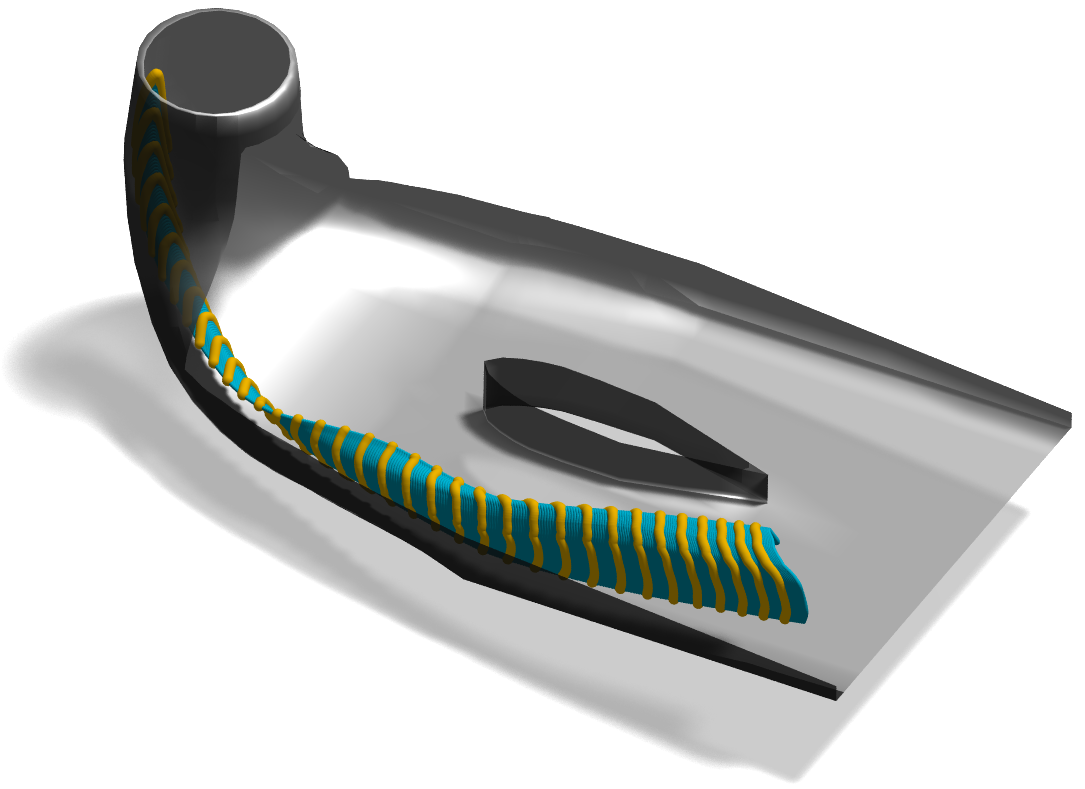}
    \put(0,5){(f)}
	\end{overpic}\hfill
\begin{overpic}[width=.24\textwidth]
    {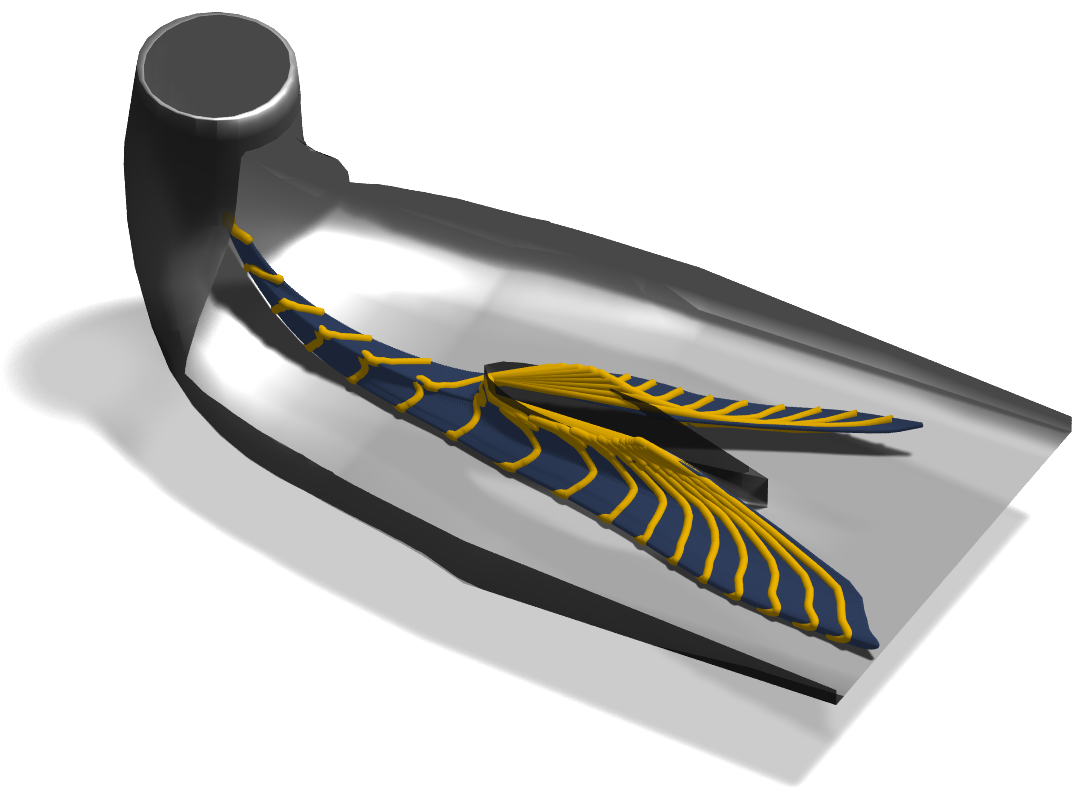}
	\put(0,5){(g)}
	\end{overpic}\hfill
\vspace{-5pt}
\end{PDF}
  \Acaption{1.5em}{Turbine data of \cite{Peikert-2002-Turbine}.
	(a) One time instant of the unsteady vector flow generated by the Francis turbine is shown.
	(b) A total of $317$ uniformly sampled points in the domain gave rise to $480$ second-order seed curves. The vector field is forward and backward integrated and the curves are ranked and colour-coded according to the strain energy $E_{\bs}$ of the associated stream surfaces; see \eqref{eq:Esrf}.
	(c) Regions with curves with low values of $E_\bs$ are further explored with higher sampling density. (d) The best-ranked initial stream surface and (e) its optimised counterpart; $E_\bs = 6.32\cdot10^{-6}$.
	(f) The second best-ranked stream surface with $E_\bs = 9.28\cdot10^{-6}$.
	(g) An example of a stream surface with a high strain energy that is split into two parts by the flow. Note that the length of the timelines grows significantly, resulting in exceeding $E_\bs = 9.28\cdot10^{2}$.}\label{fig:Turbine}
  \end{figure*}


The simulation of a flow coming from a Francis turbine is shown in Fig.~\ref{fig:Turbine}. The original vector field is unsteady, so we again used only one time instant. Note that the best strain-minimising stream surfaces are generated by seed curves that move in only one part of the turbine's body, whilst the low ranked surfaces are those where the middle part of the body forces their timelines to split apart; see Fig.~\ref{fig:Turbine}(g). This observation suggests that $E_\bs$ could be used to detect this type of splitting in a flow.

Fig.~\ref{fig:Rayleigh} shows a time instant of Rayleigh-B\'{e}nard heat convection, where the fluid is heated at the bottom and cooled at the top of the boundary, resulting in a circulatory motion. This vector field is not incompressible. Nevertheless, as seen in Fig.~\ref{fig:Rayleigh}(a),
the second-order curves still exist in this field. Recall that the guarantee of existence of first-order curves applies to divergence-free vector fields only (Remark~\ref{rem:singular}), while second-order curves are not guaranteed to exist.

\begin{figure*}[!tbh]
\begin{PDF}
 \hfill
 \begin{overpic}[width=.19\textwidth]
	{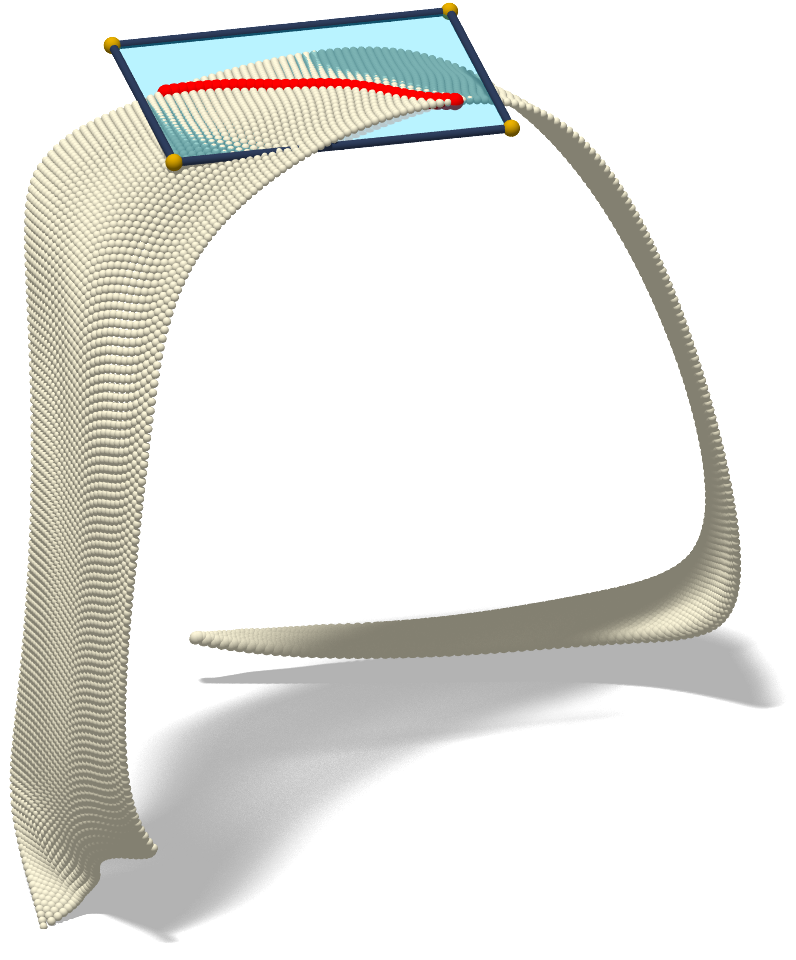}
	\put(55,95){$\partial\Omega$}
	\end{overpic}\hfill
  \begin{overpic}[width=.19\textwidth]
	{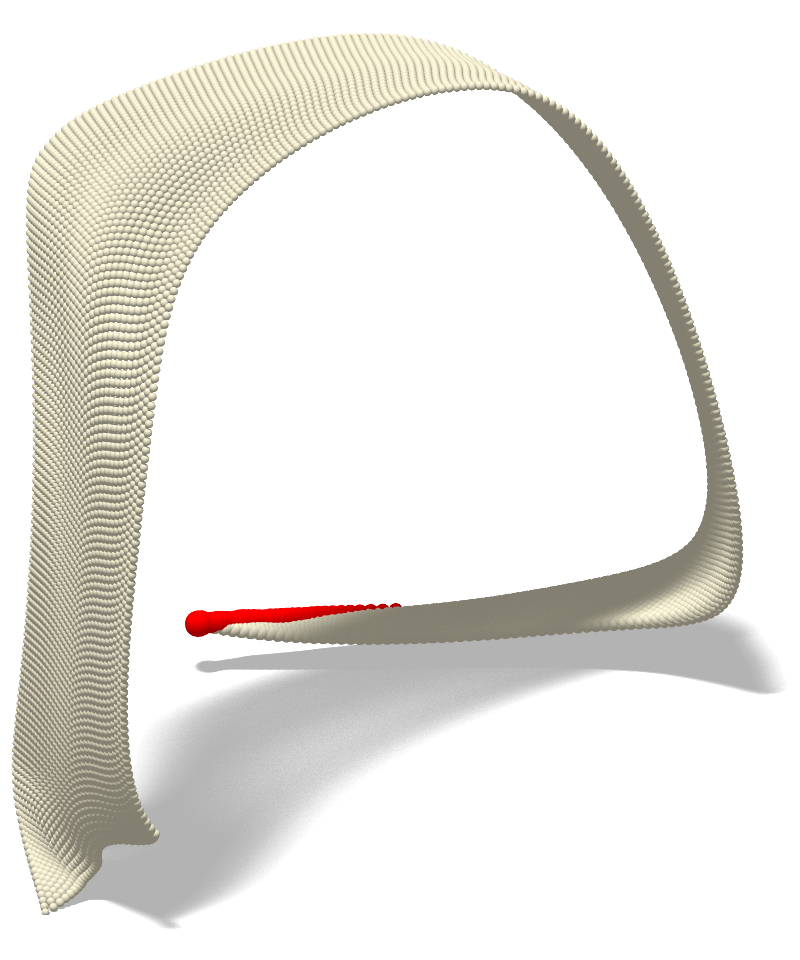}
    \put(-60,5){\small{$E^{srf} = 3.78 \cdot 10^{-4}$}}
	\end{overpic}\hfill
 \begin{overpic}[width=.19\textwidth]
	{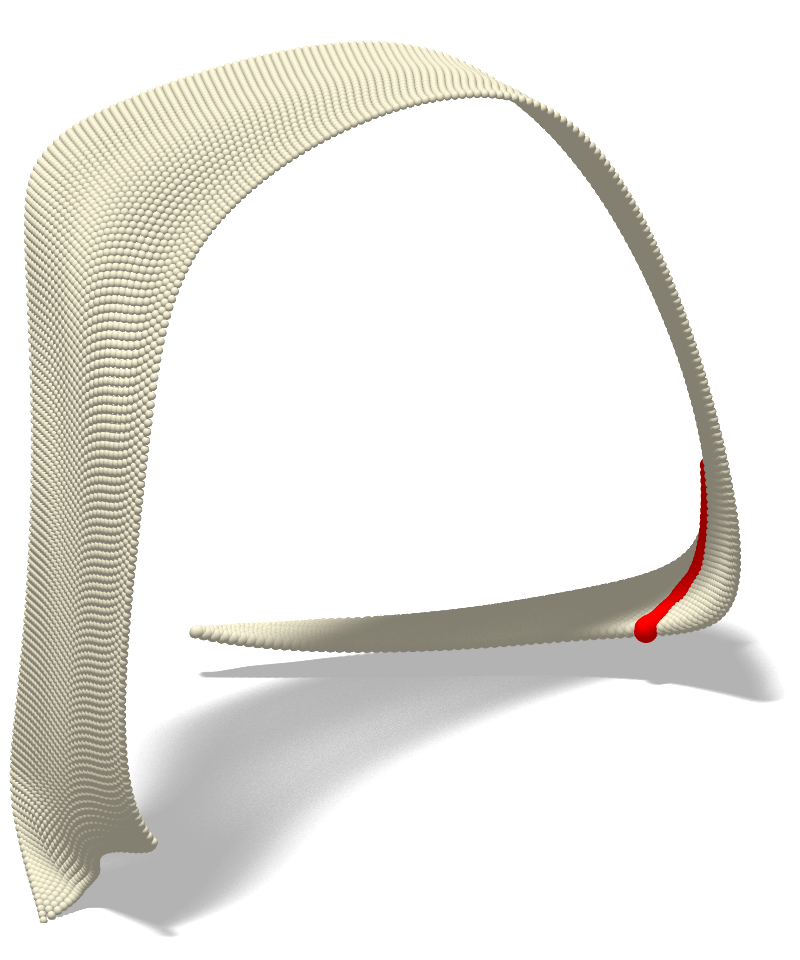}
	\end{overpic}\hfill
\begin{overpic}[width=.19\textwidth]
	{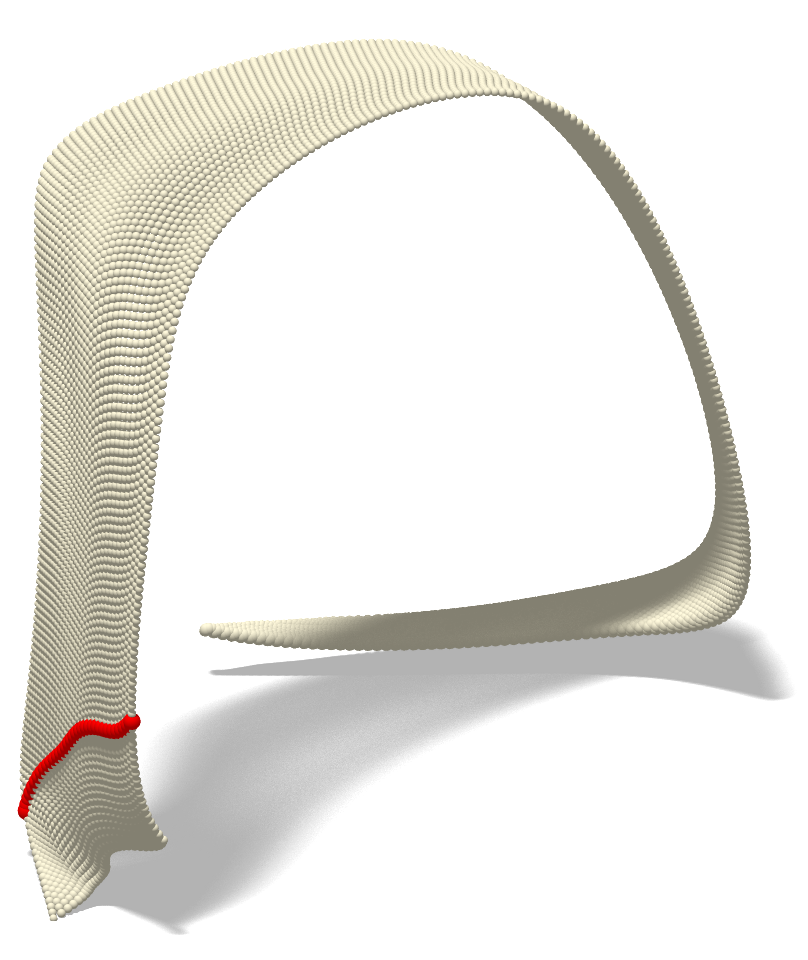}
	\end{overpic}\hfill
\begin{overpic}[width=.19\textwidth]
	{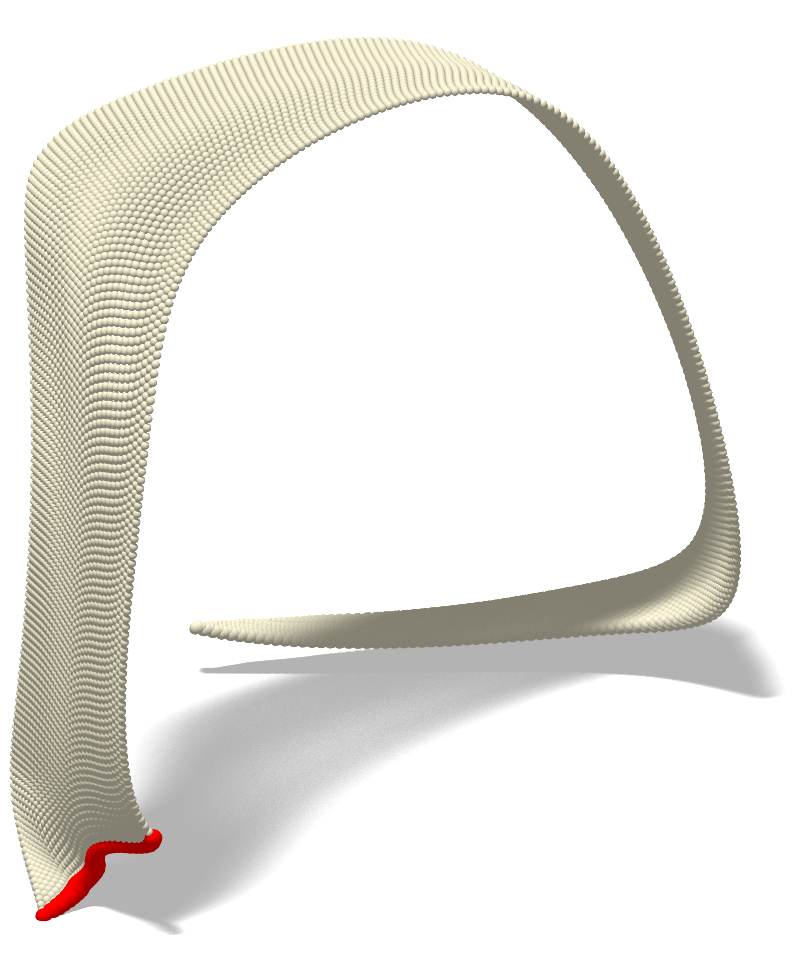}
    \put(30,-5){\includegraphics[width=0.11\textwidth]{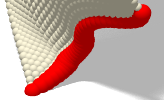}}
	\end{overpic}\hfill
\vspace{-5pt}
 \begin{overpic}[width=.19\textwidth]
	{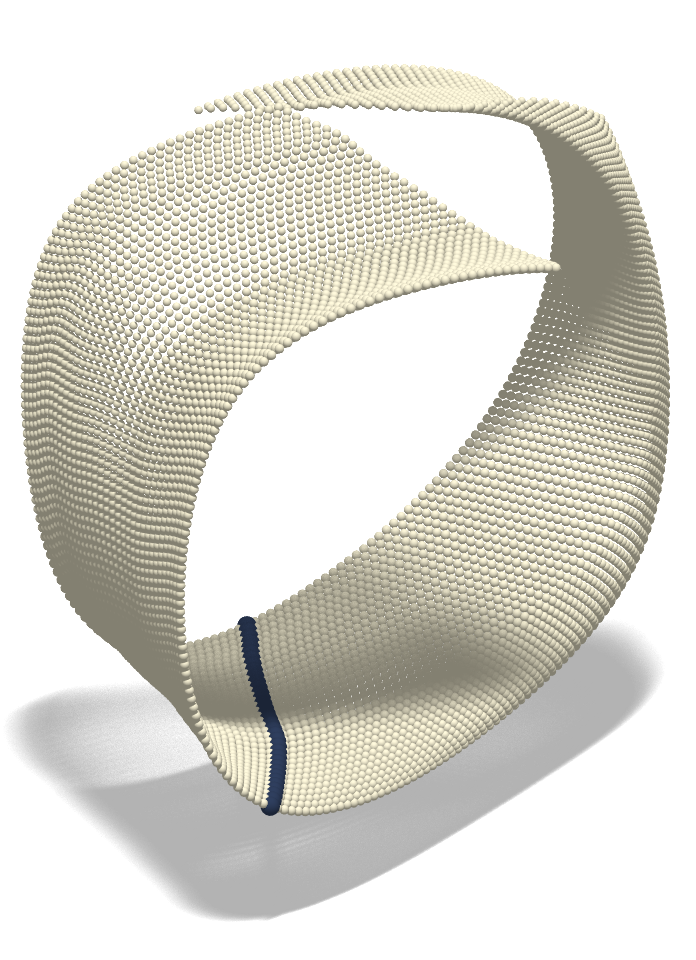}
	\put(0,90){(a)}
    \put(23,38){$\gamma$}
	\end{overpic}\hfill
  \begin{overpic}[width=.19\textwidth]
	{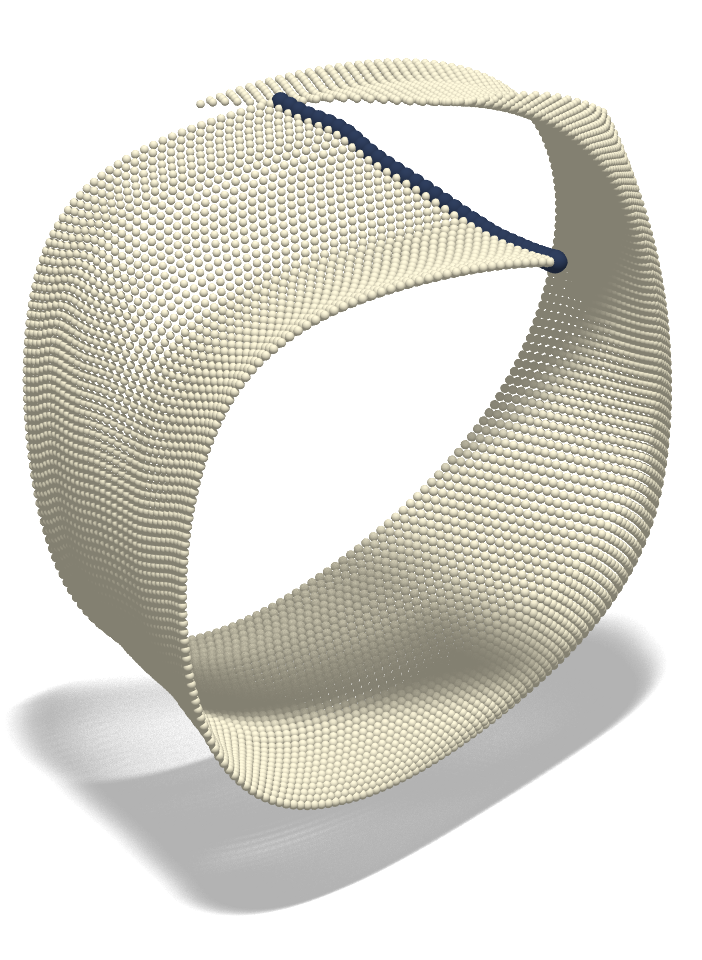}
    \put(-40,5){\small{$E^{srf} = 1.24 \cdot 10^{-5}$}}
	\put(0,90){(b)}
	\end{overpic}\hfill
 \begin{overpic}[width=.19\textwidth]
	{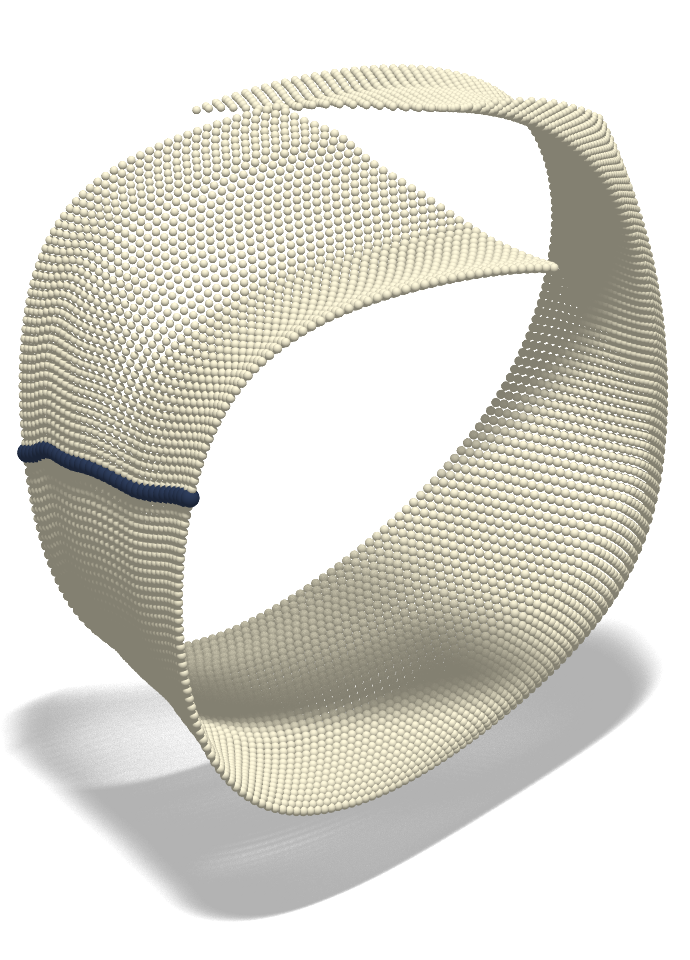}
    \put(0,90){(c)}
	\end{overpic}\hfill
\begin{overpic}[width=.19\textwidth]
	{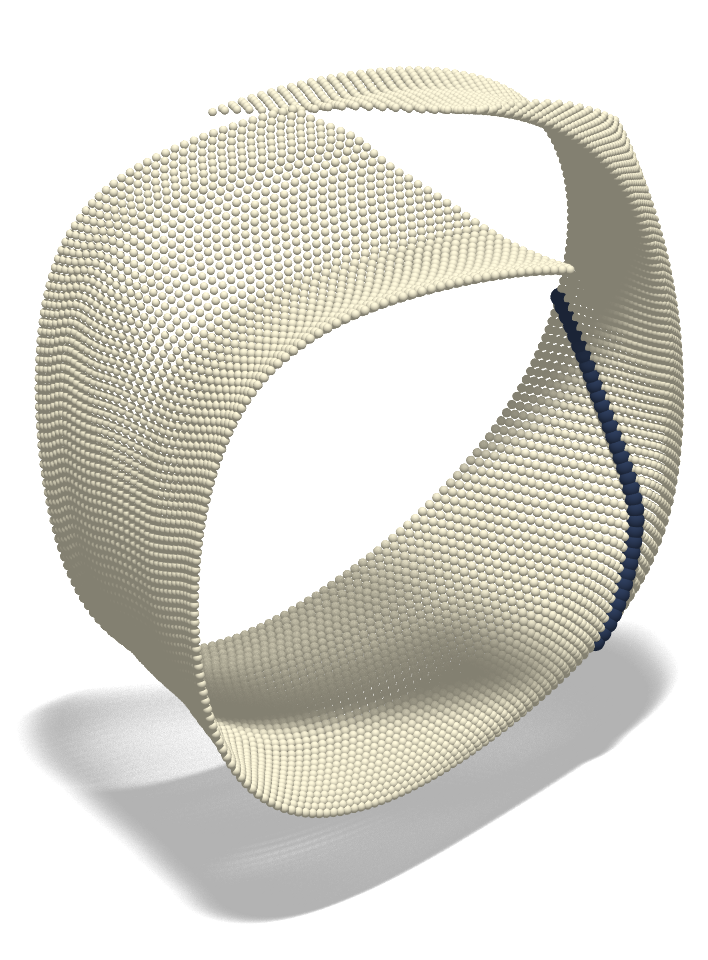}
	\put(0,90){(d)}
	\end{overpic}\hfill
\begin{overpic}[width=.19\textwidth]
	{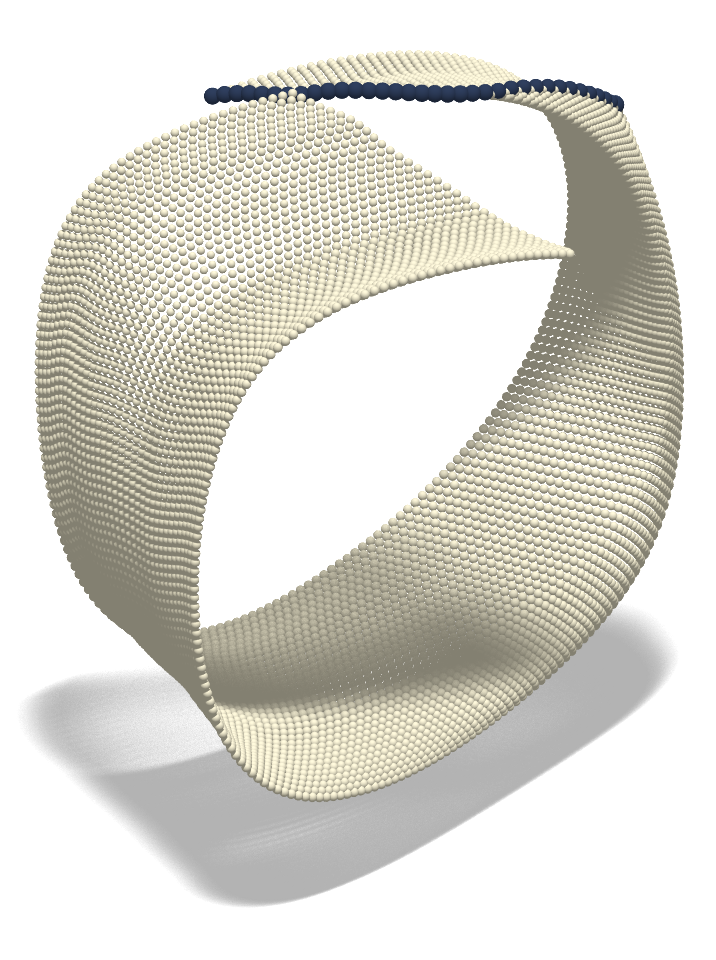}
	\put(0,90){(e)}
	\end{overpic}\hfill
\vspace{-4pt}
\end{PDF}
  \Acaption{1.5em}{A comparison of first- and second-order curves. The best initial stream surface according to $E_\bs$ given by first-order (top)
  and second-order (bottom) curves. (a) The initial first-order curve is constrained to lie on $\partial \Omega$ (top), whereas the second-order curve $\gamma$ lies inside $\Omega$ (bottom). (b--e) Visualisation of the strain energy $E_\bs$ along the stream surfaces is facilitated by spheres at uniformly sampled points along the seed curve. The spheres' radii reflect the change of the magnitude of the tangent vectors of the time curves. Observe that second-order curves serve as perfect initial guesses for optimisation. All the spheres are of almost equal radii, which corresponds to $E_\bs$ being close to zero, whilst the radius varies along the first order curve as is shown in the close-up of top (e). See also the accompanying video.}\label{fig:1stVs2nd}
  \end{figure*}

Fig.~\ref{fig:1stVs2nd} shows a comparison between first- and the second-order curves when used as seed curves in the vector field from Fig.~\ref{fig:Rayleigh}. As expected, second-order curves offer more favourable results and give rise to initial stream surfaces with nearly vanishing strain
energy $E_\bs$.

Statistics concerning the number of sampled points and seed curves, optimisation parameters, resulting energies $E_\bs$, and computation times are listed in Table~\ref{tab}.
The timings differ depending on concrete parameter settings, ranging from a few seconds to several minutes. The most expensive part is the exploration of families of first- and/or second-order curves, with computing the energy $E_{\bs}$ in particular, which requires integration of all the candidate seed curves.
The reported timings were obtained on a machine running Windows with a double-core CPU (2.67 GHz) and 24GB RAM. Currently, only a single-core CPU implementation is available. However, the algorithm is well suited for parallelisation (e.g. curve/surface energies can be computed independently) on the CPU or even GPU.

  \begin{table}[t]
 \begin{center}
  \begin{minipage}{0.99\columnwidth}
\caption{A summary of our results. In the \emph{`Seed curves'} columns, $\# \bp$ is the number of sampled points, $\# \gamma$ is the number of (first- or second-order) seed curves emanating from them, and $E_{\bs}^{ini}$ is the strain energy \eqref{eq:Esrf} of the best candidate before optimisation. In the \emph{`Optimisation'} columns, $\# \bs$ is
the number of surfaces processed and $E_{\bs}^{opt}$ is the energy of the best optimised stream surface. The last column lists total computation times in seconds.}\label{tab}
  \end{minipage}
\vspace{3pt}\\
{\footnotesize
\renewcommand{\arraystretch}{1.1}
\begin{tabular}{|c || r| r| r| r| r| r| r|}\hline
 & \multicolumn{3}{c|}{Seed curves}
 & \multicolumn{2}{c|}{Optimisation} & Time\\
 Fig. & $\# \bp$ & $\# \gamma$ & $E_{\bs}^{ini}$ & $\# \bs$ & $E_{\bs}^{opt}$ & (sec.)
\\\hline\hline
 \ref{fig:EnergyRanking}(a--c) & 30 & 41 & $6.04e^{-2}$ & n/a & n/a & 2 \\\hline
 \ref{fig:Turbine}             & 240 & 492 & $3.17e^{-4}$ & 54 & $6.32e^{-6}$ & 385  \\\hline
 \ref{fig:Rayleigh}            & 216 & 434 & $1.24e^{-5}$ & 76 & $2.44e^{-7}$ & 513  \\\hline
 \ref{fig:SquareCylinder}(b)   & 216 & 239 & $6.24e^{-4}$ & 32 & $8.17e^{-6}$ & 147  \\\hline
 \ref{fig:SquareCylinder}(c)   & 216 & 184 & $9.51e^{-4}$ & 21 & $1.04e^{-5}$ & 135  \\\hline
\end{tabular}
} \end{center} \end{table}

\section{Discussion, limitations and future work}\label{sec:dis}

We now discuss extensions, limitations and avenues for future research.

\textbf{Rigid body flow.} In our implementation, we have considered only generic cases, i.e., when $\textnormal{det}(\bj)\neq 0$ and $\textnormal{det}(\bk)\neq 0$.
 If a singular case was detected while integrating a seed curve, the integration was terminated. In the special case of rigid body flows, every curve is strain free and the problem becomes ill-posed. On the other hand, the singular cases can be easily detected.

\textbf{Unstructured grids.} All the vector fields tested in this paper were known at vertices of very fine structured grids (or analytically). This brings  certain simplifications, e.g. when estimating the vector field outside the grid-points. Our implementation, if needed, could be easily extended to accommodate unstructured volumetric meshes as well.

\textbf{Fields with non-vanishing divergence.} We have focused on divergence-free vector fields since this guarantees the existence of strain minimising curves (of first order). Nevertheless, as shown in Fig.~\ref{fig:Rayleigh}, such seed curves may exist even in general vector fields and explorations in this direction seem promising.

\textbf{Unsteady vector flows.} For the sake of simplicity, only steady vector fields were considered. The generalisation to the unsteady case is
straightforward. This would extend the dimensionality of the space of seed curves from three to four.

\textbf{Non-local first-order curves.} One could consider `non-infinitesimally' arc-length minimising/vanishing seed curves, i.e., curves that preserve their arc-length after a \emph{non-zero} time increment $\Delta t$.
The problem of finding the intersection of two quadratic cones, as encountered in Section~\ref{sec:SE2}, would turn into a sphere-ellipsoid intersection problem (induced by the linear map between the tangent spaces at $t=0$ and $t=\Delta t$), which is more demanding. Also, it is not obvious how to set an appropriate $\Delta t$.

\textbf{Area-preserving surfaces.} A natural generalisation of our method points towards area-preserving surfaces, i.e., surfaces that propagate in time
while preserving the area of any of their sub-patches. As the family of first-order curves is large, it is reasonable to expect that such surfaces exist.

\textbf{Stability.} As a stream surface is determined by its seed curve, the seeding strategy is a crucial ingredient. We have tested random and uniform seeding schemes. Whereas the uniform sampling was applied in most of the examples, followed by adaptive sampling in the neighborhoods with low energy values, in Fig.~\ref{fig:EnergyRanking} the test with random sampling of the boundary points and consequently boundary seeding curves is shown to validate the stability of the algorithm.

\section{Conclusion}\label{sec:conclu}

Combining theoretical and numerical techniques, we have attacked the difficult problem (with no known close-form solution) of finding curves that propagate in an arc-length preserving manner in divergence-free vector fields and thus give rise to strain-minimising stream surfaces. Our method successfully solves this problem.

Our contributions are both theoretical and practical. On the one hand, we have advanced the theory of strain-minimising curves and surfaces in divergence-free vector fields. We have shown that such objects always exist. On the other hand, we have illustrated by several examples and the accompanying video that our theoretical results lead to immediate applications such as vector field exploration and visualisation.

\paragraph*{Acknowledgments}
We would like to thank Mike Schulze for providing the Turbine \cite{Peikert-2002-Turbine} and the `Square cylinder' \cite{Weinkauf-2008-SquareCylinder} datasets, and Holger Theisel for his suggestions and comments.
The second author were supported by EPSRC through Grant EP/H030115/1.

\section*{Apendix A}

We now prove Lemma~3.1. By definition,
$ c_1 = C_1(0)$, where
$C_1(t) = \frac{\partial }{\partial t}\int_{s_0}^{s_1} ||\bs_s(s,t)|| \, \mathrm{d}s$. We have that
\begin{equation*}
\begin{array}{rcl}
C_1(t) &=& \int_{s_0}^{s_1} \frac{\partial }{\partial t}||\bs_s(s,t)|| \, \mathrm{d}s = \int_{s_0}^{s_1} \frac{\partial }{\partial t}\sqrt{\dpr{\bs_s}{\bs_s}} \, \mathrm{d}s\\
&=& \int_{s_0}^{s_1} \frac{\dpr{\bs_s}{\bs_{st}}}{\sqrt{\dpr{\bs_s}{\bs_s}}} \, \mathrm{d}s.
\end{array}
\end{equation*}
Consequently, since $\sqrt{\dpr{\bs_s}{\bs_s}}|_{t=0}\equiv 1$, we can conclude that
$c_1 = \int_{s_0}^{s_1}\dpr{\bs_s}{\bs_{st}} \, \mathrm{d}s |_{t=0}$.
Proceeding similarly in the case of $c_2$, we have
\begin{equation*}
\begin{array}{rcl}
C_2(t) &=& \frac{\partial^2 }{\partial t^2}\int_{s_0}^{s_1} ||\bs_s(s,t)|| \, \mathrm{d}s =
\int_{s_0}^{s_1} \frac{\partial^2 }{\partial t^2}\sqrt{\dpr{\bs_s}{\bs_s}} \, \mathrm{d}s\\
&=& \int_{s_0}^{s_1} \frac{\dpr{\bs_{st}}{\bs_{st}}+\dpr{\bs_{s}}{\bs_{stt}}}
                        {\sqrt{\dpr{\bs_s}{\bs_s}}}
                 - \frac{\dpr{\bs_{s}}{\bs_{st}}^2}
                        {\sqrt{\dpr{\bs_s}{\bs_s}}^3} \, \mathrm{d}s.
\end{array}
\end{equation*}
Therefore,
$$c_2 = C_2(0) = \int_{s_0}^{s_1}\dpr{\bs_{st}}{\bs_{st}}+\dpr{\bs_s}{\bs_{stt}}-\dpr{\bs_s}{\bs_{st}}^2 \, \mathrm{d}s |_{t=0}$$
as claimed. \hfill $\square$

\section*{Appendix B}

We now present a geometric approach to the problem of computing the intersection vectors of two quadratic cones given by $\bd \bj \bd^\top = 0$ and $\bd \bk \bd^\top = 0$ with apexes at $\bp$. The two cones belong to a pencil of cones given by $\bd(\bj+\lambda\bk)\bd^\top=0$ parametrised by $\lambda$. We identify the singular quadric in this pencil by setting $\mathrm{det}(\bj+\lambda\bk)=0$, which, generically, leads to a cubic equation in $\lambda$. Thus, there exists at least one real $\lambda_0$ which determines a singular quadric in the family. If the intersection of the two cones is real and non-trivial (i.e., not $\bp$ or a cone), the real part of the singular quadric given by $\bj+\lambda_0\bk$ is a pair of planes (possibly coincident) or a straight line, and is incident with $\bp$. Finally, intersecting this line or the planes with either of the input cones is a simple quadratic problem which gives the sought-after second-order vector(s), provided that they exist.

\section*{Appendix C}
\begin{figure}[!tbh]
 \hfill
  \begin{overpic}[width=.95\columnwidth]
	{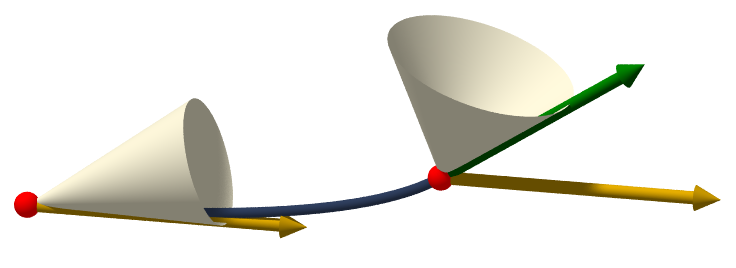}
	\put(-2,13){$\bp(t_0)$}
    \put(40,2){$\bd(t_0)$}
    \put(58,6){$\bp(t_1)$}
    \put(88,25){$\bd(t_1)$}
    \put(20,25){$\bd \bj(t_0) \bd^\top=0$}
    \put(67,33){$\bd \bj(t_1) \bd^\top=0$}
	\end{overpic}\hfill
    \vspace{-2pt}
  \Acaption{1.5em}{Definition and computation of first-order interior curves. A first-order vector $\bd$ at $t=t_0$ (yellow) is projected (green) onto the first-order cone at $t_1=t_0+h$.}\label{fig:StraightCone}
  \end{figure}

The integration procedure for first-order interior curves is shown in Fig.~\ref{fig:StraightCone}. Denote $h$ the step-size (set to $\mathrm{diam}(\Omega)/10$ by default). Given a point $\bp(t_0)\in\Omega$ and an initial first-order vector $\bd(t_0)$, the first-order vector at $t_1=t_0+h$ is obtained by projecting $\bd(t_0)$ translated to $\bp(t_1)$ onto the quadratic cone given by $\bd \bj(t_1) \bd^\top=0$.

\bibliographystyle{eg-alpha-doi}
\bibliography{StreamSrfArXiv}

\end{document}